\documentclass[nonblindrev]{informs3b}
\SingleSpacedXI
\TheoremsNumberedThrough
\EquationsNumberedThrough
\usepackage[T1]{fontenc}
\usepackage{mathtools,booktabs,subcaption,float}
\usepackage[margin=1in]{geometry}
\usepackage{natbib}
\bibpunct[, ]{(}{)}{,}{a}{}{,}
\def\bibfont{\small}

\setcitestyle{authoryear,open={(},close={)}}
\usepackage[hidelinks]{hyperref}
\hypersetup{
  pdftitle={When Suspicion Adds Up: A Stochastic Framework for Optimizing Honeypot-Based Insider-Threat Detection},
  pdfauthor={Saharnaz Mehrani, Lujia Zhan, Mohsen Emadikhiav, Payman Dehghanian},
  pdfsubject={Insider threat detection with honeypots},
  pdfkeywords={Insider threat detection; Cyber-physical systems; Continuous-time Markov chain; Stochastic optimization.}
}
\renewcommand{\ABSfont}{\small}
\renewcommand{\AFFfont}{\small}
\begin{document}
\raggedbottom
\RUNAUTHOR{Mehrani et al.}
    \RUNTITLE{A Stochastic Framework for Insider Threat Detection}
	\TITLE{When Suspicion Adds Up: A Stochastic Framework for Optimizing Honeypot-Based Insider-Threat Detection}
	\ARTICLEAUTHORS{
\AUTHOR{Saharnaz Mehrani\textsuperscript{1}, Lujia Zhan\textsuperscript{2},
Mohsen Emadikhiav\textsuperscript{1}, Payman Dehghanian\textsuperscript{2}}
\AFF{\textsuperscript{1}Department of Information Technology and Operations Management,
College of Business, Florida Atlantic University,\\
777 Glades Road, Boca Raton, FL 33431, USA}
\AFF{\textsuperscript{2}Department of Electrical and Computer Engineering,
School of Engineering and Applied Science, The George Washington University,\\
800 22nd Street NW, Washington, DC 20052, USA}
\AFF{\href{mailto:smehrani@fau.edu}{smehrani@fau.edu},
\href{mailto:lujiacky@gwmail.gwu.edu}{lujiacky@gwmail.gwu.edu},
\href{mailto:memadikhiav@fau.edu}{memadikhiav@fau.edu},
\href{mailto:payman@gwu.edu}{payman@gwu.edu}\\
Corresponding author: Saharnaz Mehrani}
}

    \ABSTRACT{
    Insider threats pose security risks when employees or other trusted users misuse their access to organizational systems. Detecting these threats is difficult because harmful activity can resemble routine behavior. Honeypots, or decoy resources monitored for suspicious interactions, can provide evidence of potentially malicious activity. We develop a decision framework that models evidence accumulation through honeypot interactions using a continuous-time Markov chain and jointly selects honeypot counts and monitoring horizons to minimize expected losses from malicious activity and false detections. For thresholds of two and three interactions, we derive closed-form true- and false-detection probabilities and expected detection times conditional on detection within the monitoring horizon. We characterize boundary cases and establish stochastic-ordering results. Our numerical results highlight the value of the proposed framework across different settings. %

    }

\maketitle 

\section{Introduction}\label{sec:Introduction}

Modern \emph{cyber--physical systems} (CPSs) integrate networked sensing, computation, communication, and control with physical operations. This integration has improved the visibility and automation of critical infrastructures, including electric power and transportation systems, but it has also created channels through which cyber compromises can affect physical processes~\citep{6016202,zhu2011taxonomy,li2017distributed}. Insider threats are especially difficult to address because employees and other trusted users can misuse legitimate access and knowledge of system operations, making harmful activity harder to distinguish from routine behavior~\citep{homoliak2019insight}. In practice, harmful intent may not be visible at the outset, and motives are difficult to observe directly. Organizations therefore rely on controls that reduce opportunities for misuse and encourage security-policy compliance~\citep{padayachee2016assessment,herath2009encouraging}. Research also suggests that insider maliciousness can develop through a progression of stages~\citep{johnston2024cyber}. These considerations create a need for detection mechanisms that monitor behavior over time and identify misuse before prolonged exposure causes substantial operational harm.

\begin{figure}[!t]
\centering
\includegraphics[width=0.6\linewidth]{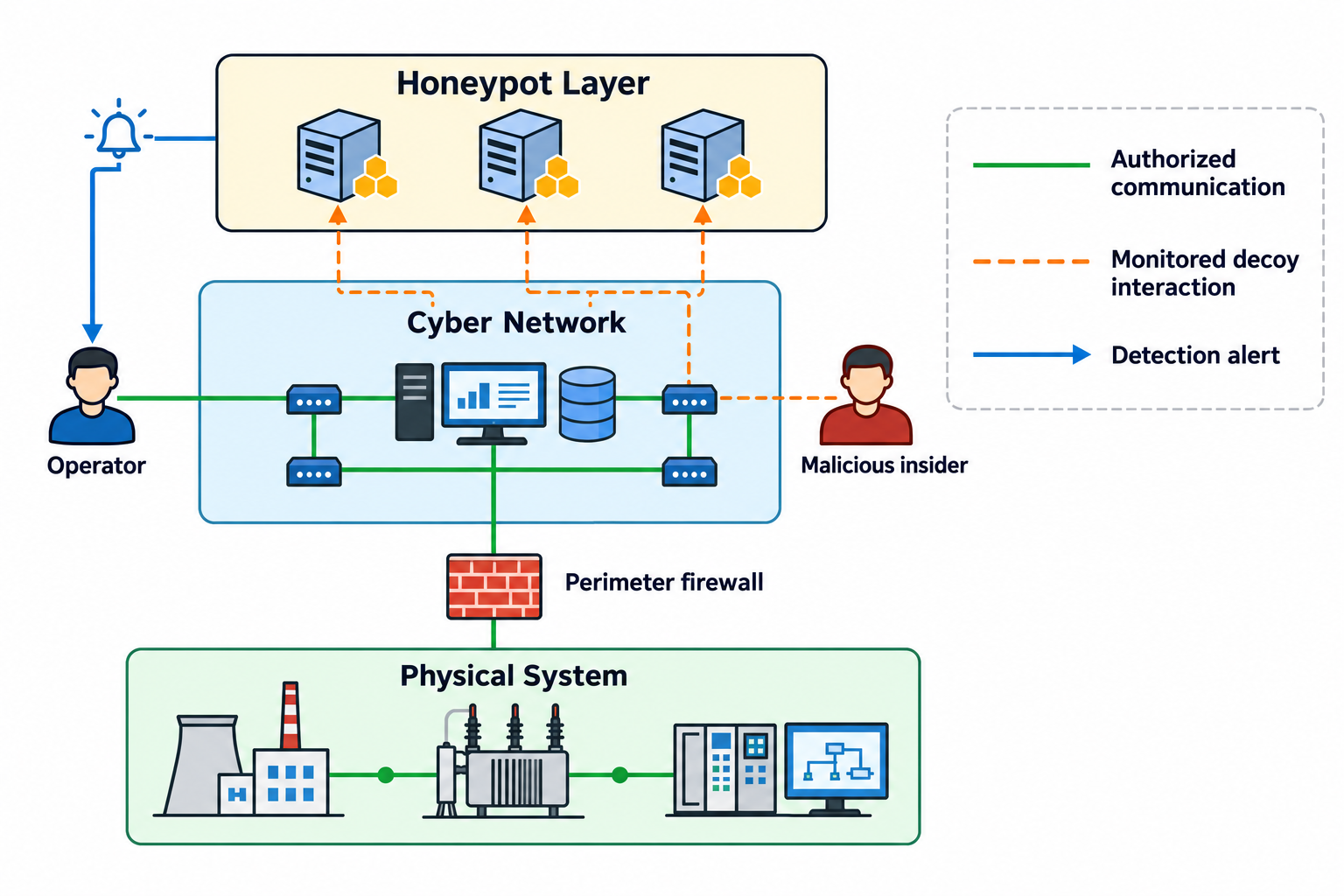}
\caption{Illustration of honeypots within a CPS.}
\label{fig:honeypot_role_cpps}
\end{figure}

Honeypots offer a complementary mechanism for observing such misuse. A honeypot is an instrumented decoy resource whose security value arises from monitored interaction~\citep{spitzner2003honeypots}. Because normal operational tasks do not require engagement with these decoys, interactions with them can provide focused evidence of suspicious behavior. Honeypots therefore complement firewalls, intrusion detection systems (IDSs), and endpoint protection. An interaction alone, however, does not establish malicious intent. A benign insider may encounter a decoy inadvertently, while a malicious insider may engage repeatedly, interact only sporadically, or recognize and avoid it. These possibilities motivate evaluating the timing and history of honeypot interactions. In a power-system setting, for example, honeypots can imitate load measurements, control interfaces, commands, status information, and actuator responses~\citep{zhan2026minimizing}. Interactions with these decoys can reveal attempts to access or manipulate system information and prompt restrictions on insider access. Figure~\ref{fig:honeypot_role_cpps} illustrates the role of honeypots as decoy assets within a CPS.

Repeated interactions create an operational trade-off involving the number of honeypots, the duration of monitoring, and the evidence required for detection. Deploying more honeypots creates additional opportunities to observe malicious activity and can accelerate detection, but it can also increase the likelihood that benign interactions accumulate to the detection threshold. A longer monitoring horizon allows more time to collect evidence and may support a smaller deployment, but it also gives benign interactions more time to reach the detection threshold. Conversely, a short horizon may end before sufficient evidence is collected. Delayed or missed malicious detection creates exposure costs, making detection timing relevant to these decisions~\citep{paul2021decision,safi2023detecting}. Because the same detection rule is applied to both insider types, these choices affect the likelihood of detecting malicious insiders, the likelihood of falsely accusing benign insiders, and the time required for detection. We therefore address the following question: how should the number of honeypots and the monitoring horizon be chosen to support timely detection of malicious insiders while controlling false detections?

We develop a stochastic framework in which detection occurs after an insider completes a specified number of honeypot interactions within a single finite monitoring horizon. Throughout the horizon, the insider's type is fixed as benign or malicious. State transitions represent the accumulation of interaction evidence rather than a change in the insider's underlying intent. We model this process through a continuous-time Markov chain with type-specific interaction rates. The rates distinguish first interactions with previously unvisited honeypots from repeated visits, allowing the model to capture differences in both interaction frequency and persistence. %

This study examines how repeated engagement with individual honeypots affects deployment and monitoring decisions. We develop this contribution in four parts. First, we distinguish first visits from repeated visits through type-specific interaction rates and evaluate benign and malicious insiders under the same detection rule. %
We examine two- and three-stage models, corresponding to detection thresholds of two and three interactions. We focus on these cases because two interactions provide the simplest threshold that requires more than one interaction, while three interactions require distinguishing whether the first two visits involved the same or different honeypots. This distinction preserves the relevant interaction history in the three-stage model. Second, for both models, we derive closed-form expressions for the true- and false-detection probabilities and the expected detection time conditional on detection within the monitoring horizon. We characterize boundary cases, establish stochastic ordering of malicious and benign detection times under cross-type rate-dominance conditions, and determine how the higher evidence threshold changes detection performance. Third, we embed these analytical measures in an optimization model that jointly selects the honeypot count and monitoring horizon for a given evidence threshold. The objective balances the exposure costs of delayed or missed malicious detection against the costs of falsely detecting benign insiders.

Fourth, our numerical analysis examines how interaction rates and relative costs affect deployment and monitoring decisions, and compares evidence thresholds after those decisions are adjusted. Additional honeypots can reduce malicious exposure even when detection probabilities are nearly saturated. The false-detection advantage of a higher threshold at a fixed configuration can also disappear when deployment and monitoring are adjusted to meet a common malicious-detection target. Under cost minimization, the preferred threshold depends on the relative costs of malicious exposure and falsely accusing benign insiders. These findings show why evidence requirements should be evaluated together with deployment size and monitoring duration.

The remainder of the paper is organized as follows. Section~\ref{sec:LiteratureReview} reviews the related literature and positions the study. Section~\ref{Sec:Framework} presents the detection framework, derives the performance measures, and formulates the optimization problem. Section~\ref{sec:TheoreticalProperties} establishes the structural properties of the two- and three-stage models. Section~\ref{sec:experiments} reports the numerical experiments and discusses the resulting detection and design trade-offs. Section~\ref{sec:conclusions} concludes the paper and outlines directions for future research.

\section{Literature Review}\label{sec:LiteratureReview}

Our study relates to research on insider behavior and threat detection, honeypot-based defense, and security monitoring and resource allocation. This section reviews the studies most relevant to our work and clarifies the gap our study addresses. Existing research establishes methods for accumulating detection evidence, adapting decoy defenses, and configuring detection systems on economic grounds. Our study focuses on how repeated engagement with individual honeypots changes deployment and monitoring decisions. We distinguish first and repeated visits by benign and malicious insiders and link detection probability and timing to the joint choice of honeypot count and monitoring horizon. This allows evidence thresholds to be compared after deployment and monitoring are adjusted, accounting for the costs of malicious exposure and false detections.

\noindent\textbf{Organizational controls and insider behavior.}
Research on security behavior examines deterrent and preventive measures against computer abuse \citep{straub1990effective}, penalties and social pressures affecting security-policy compliance intentions \citep{herath2009encouraging}, and opportunity-reducing security techniques \citep{padayachee2016assessment}. \citet{okeke2023application} identify unclear organizational roles and shared data-security responsibilities as barriers to preventing internal identity theft in UK retail firms. Other studies examine insiders' perceptions of risks and benefits \citep{farahmand2013understanding}, motives for computer abuse \citep{burns2023going}, and progression toward malicious behavior \citep{johnston2024cyber}. These studies highlight the behavioral context of insider threats and the importance of organizational safeguards.

\noindent\textbf{Behavioral monitoring and evidence accumulation.}
Insider threat detection draws on technical observations and behavioral and organizational information \citep{homoliak2019insight}. \citet{greitzer2010combining} motivate combining security audit data with psychosocial indicators, while \citet{legg2015automated} use individual and role-based profiles to identify anomalous behavior. Learning-based approaches include online anomaly detection with interpretable feature contributions \citep{tuor2017deep} and behavioral deep learning \citep{nasir2021behavioral}. In related work on network anomalies, \citet{padin2026toward} generate natural-language explanations of alerts to support analysts' interpretation. Some studies explicitly focus on evidence accumulation over time. \citet{chivers2013knowing} use incremental Bayesian updates to identify potential attackers whose activities are obscured by false alarms and background noise. \citet{panigrahi2013two} combine evidence within and across database transactions through belief updating and evaluate detection performance using stochastic models. Our focus is how first and repeated visits to individual honeypots affect detection timing and operational decisions.

\noindent\textbf{Honeypot-based detection and defense.}
Honeypots elicit and reveal suspicious activity through decoy assets. \citet{spitzner2003honeypots} discusses their use against insiders through honeytokens and monitored decoy environments, while \citet{yamin2019implementation} propose encrypted honeypot sensors to support insider detection and reduce false positives. Analytical studies consider defense architecture and adaptation. \citet{shi2018performance} use stochastic Petri nets and associated Markov chains to compare defense configurations. \citet{shi2019dynamic} develop a blockchain-based distributed architecture that changes real and decoy service roles across hosts. \citet{shi2021research} analyze array-honeypot defense through a three-party evolutionary game, and \citet{wang2024game} combine evolutionary game dynamics and learning to adapt real and decoy service switching. Closely related to our setting, \citet{zhan2026minimizing} jointly select component protection and honeypot counts to reduce the financial risk of load-redistribution attacks involving insiders. Their model uses a specified per-honeypot interaction probability and assumes that one interaction triggers insider removal. We examine how waiting for additional interactions changes detection timing and false detections, and consequently the preferred honeypot count and monitoring horizon.

\noindent\textbf{Security monitoring and resource decisions.}
Economic models connect detection-system settings to their operational costs. \citet{ccakanyildirim2009management} develop multi-period models that inform IDS configuration, alarm investigation, and investment while accounting for investigation costs, ignored alarms, and undetected intrusions. Complementing these models, \citet{safi2023detecting} experimentally show that recent incidents alter threat-checking patterns and that added protection can weaken monitoring, although this response does not fully offset its benefits. \citet{wu2026fixed} examine security-service contracts, showing how contract terms shape client negligence and affect providers' market coverage and profitability.

Broader planning models address control selection and investment timing. \citet{zhang2018decision} allocate security controls under uncertain breach probabilities while accounting for business impacts. \citet{khouzani2019scalable} optimize controls against multistage attacks represented by probabilistic attack graphs. \citet{zhang2021bayesian} combine preventive portfolio optimization and learning about ongoing attacks with online optimization using Bayesian Stackelberg games. \citet{benaroch2018real} uses real options to evaluate and sequence uncertainty-reducing mitigations, while \citet{paul2021decision} account for detection delays when coordinating firms' security investments with government intelligence investments. \citet{riosinsua2021adversarial} provide a broader risk-analysis framework linking uncertain threats and impacts to security-control and insurance decisions, with risk preferences and strategic attacker behavior. Our model could inform the evaluation of honeypot detection within such a process. %

\noindent\textbf{Summary of our contributions.}
Our contribution lies in connecting repeated engagement with individual honeypots to deployment and monitoring decisions. The model distinguishes first and repeated visits through type-specific interaction rates and evaluates benign and malicious insiders under the same detection rule. For thresholds of two and three interactions, we derive closed-form expressions for true- and false-detection probabilities and expected detection times conditional on detection within a single monitoring horizon. Boundary cases and stochastic-ordering results characterize the resulting detection processes. These measures support joint selection of the honeypot count and monitoring horizon to balance malicious exposure costs from delayed or missed detection against false-detection costs. After optimizing deployment and monitoring separately for each threshold, we compare the resulting costs with a first-interaction benchmark. The numerical results show that waiting for more evidence can lower modeled cost, while high exposure costs relative to false-detection costs can favor detection at the first interaction.

\section{Proposed Detection Framework}\label{Sec:Framework}

To effectively model insider threat detection in CPSs, it is essential to capture the temporal evolution of insider behavior. Insider actions are inherently dynamic, as they evolve while insiders continuously interact with system components. Observing this behavioral evolution is critical for distinguishing benign activity from emerging malicious intent and for enabling timely detection of insider-driven security breaches. In particular, a malicious insider may exhibit increasing activity levels, context-dependent deviations, or anomalous interaction patterns in response to new, sensitive, or unexpected system information, whereas a benign insider's behavior typically remains stable, predictable, and consistent over time.

To capture this dynamic process, our detection framework models the accumulation of evidence from honeypot interactions over a finite monitoring horizon. An insider's underlying type remains fixed as benign or malicious within the horizon, while the system's detection status changes as interactions are observed. Type-specific interaction rates represent differences in the frequency of first visits and repeated engagement with honeypots.

\begin{figure}[t]
\centering
\setlength{\abovecaptionskip}{-1mm}
\includegraphics[width=1\columnwidth]{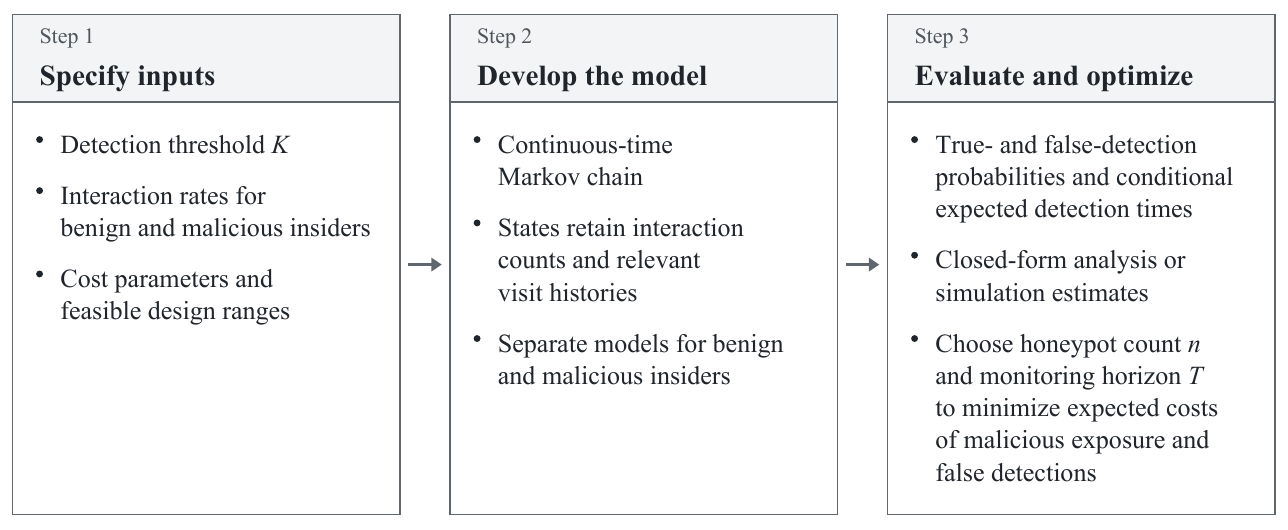}
\caption{Steps within our insider threat detection framework}
\vspace{-1mm}
\label{fig: Detection framework}
\end{figure}

The framework consists of three main steps, as illustrated in Fig.~\ref{fig: Detection framework}. In the first step, the system operator specifies (i) the number of honeypot interactions required for an insider to be flagged as malicious, denoted by $K\ge 1$, and (ii) the interaction-rate parameters for benign and malicious insiders. The threshold $K$ determines both the detection capability of the honeypots in distinguishing malicious from benign insiders and the desired level of conservatism in protecting the system against insider threats. The rate parameters are specified separately for malicious and benign insiders and are influenced by the operational context as well as the informational content embedded within each honeypot (e.g., access credentials, control commands, or sensor data). This probabilistic characterization forms the foundation for quantifying behavioral differences across insider types.

In the second step, a Markov chain represents the accumulation of honeypot interactions and the corresponding detection status. Its states record the number of interactions and, where needed, the visit history that determines subsequent interaction rates. This structure permits analytical evaluation of detection probabilities and detection times for each insider type.%

In the third step, the framework evaluates detection performance for alternative honeypot counts and monitoring horizons. The main analysis focuses on $K=2$ and $K=3$, which capture evidence beyond the first interaction while allowing closed-form analysis of repeated visits. For these models, we derive analytical expressions for true- and false-detection probabilities and the expected detection time conditional on detection within the horizon. Appendix~\ref{app:K1_benchmark} provides the $K=1$ derivation for the benchmark comparison in Section~\ref{sec:K1_benchmark}. For larger values of $K$, the additional states needed to represent interaction histories may make closed-form analysis impractical. In such cases, simulation can estimate the same performance measures from repeated realizations of the detection process. These measures enter an optimization model that selects the honeypot count and monitoring horizon to minimize the expected costs of malicious exposure and false detections. Through this process, the framework provides a quantitative basis for designing cost-efficient insider threat detection strategies in CPS environments.

\subsection{Continuous-time Markov Chain}\label{Sec:CTMC}

We model the insider detection process using a continuous-time Markov chain (CTMC). Continuous time accommodates honeypot interactions occurring at irregular intervals, while discrete states track the accumulation of detection evidence. Each state retains the visit history that affects subsequent interaction rates, as illustrated by the three-stage model below. Conditional on insider type and this recorded history, we assume that transition rates remain constant during the monitoring horizon. A state with a positive total departure rate therefore has an exponentially distributed holding time. This assumption treats elapsed time in the current state as providing no additional information about the next interaction beyond that captured by insider type and the recorded history. Further, this formulation allows rates to differ between benign and malicious insiders and between first and repeated visits, while supporting tractable analysis of detection probabilities and expected detection times conditional on detection within the monitoring horizon.

The detection process is evaluated over a monitoring horizon of length $T>0$. A longer horizon increases the opportunity to accumulate evidence for both insider types, while a shorter horizon may end before the detection threshold is reached. Detection occurs as soon as the threshold is reached, without waiting for the horizon to end. Time and rates use consistent units. For example, $T=4$ represents a four-week detection horizon when rates are measured per week, with interactions occurring in continuous time. In applications, interaction rates can be estimated from monitored interactions. The exponential holding-time and independence assumptions should also be assessed against observed behavior.

\begin{figure}[t]
\centering
\setlength{\abovecaptionskip}{-1mm}
\includegraphics[width=0.75\columnwidth]{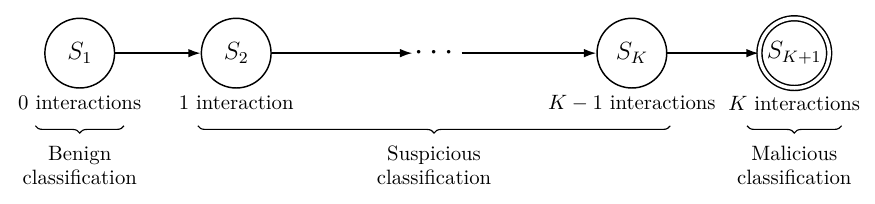}
\caption{Multi-stage model for insider-threat detection.}
\vspace{-1mm}
\label{fig: Markov Chain}
\end{figure}

Let $\mathcal{S}=\{S_1,S_2,\ldots,S_{K+1}\}$ denote the detection levels illustrated in Fig.~\ref{fig: Markov Chain}. At the start of the horizon, every insider is assigned to $S_1$, labeled \emph{Benign}, with no recorded honeypot interactions. This label denotes the initial detection status, not the insider's true type. For $K\ge2$, the first interaction takes the insider to $S_2$, labeled \emph{Suspicious}. More generally, level $S_k$ records $k-1$ interactions for $k\in\{2,\ldots,K\}$. The $K$th interaction triggers classification as \emph{Malicious} and entry into the absorbing detection state $S_{K+1}$. Intermediate levels are refined into substates when different visit histories imply different future rates. We define separate CTMCs for benign and malicious insiders using the same state representation and detection rule but type-specific rate parameters.

The analysis concerns detection within a single monitoring horizon. An insider who has not reached the detection threshold by $T$ remains undetected. This outcome does not establish that the insider is benign or imply a change in the insider's underlying type.

For insider type $(\cdot)\in\{b,m\}$, let $\tau^{(\cdot)}$ denote the first hitting time of the detection state $S_{K+1}$ for the CTMC starting from $S_1$ at time $0$. Set $\tau^{(\cdot)}=\infty$ if the threshold is never reached. The principal performance measures are $\mathbb{P}(\tau^{(\cdot)}\le T)$ and $\mathbb{E}[\tau^{(\cdot)}\mid\tau^{(\cdot)}\le T]$, with the conditional expectation defined only when detection has positive probability. In the following sections, we derive these measures and the transient state probabilities for $K=2$ and $K=3$, referred to as the two-stage and three-stage models, respectively.

\subsubsection{Two-stage Model}\label{Sec:Two-Step CTMC}

\begin{figure}[ht]
  \centering
  \begin{subfigure}[t]{0.4\linewidth}
\centering
\includegraphics[width=\linewidth]{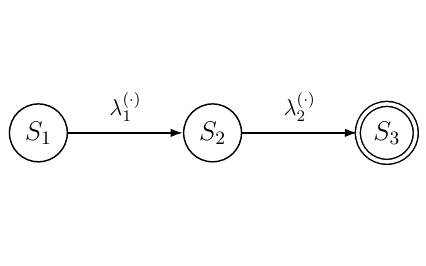}
\caption{$K=2$}
\label{fig:CTMCModelK2}
\end{subfigure}
      \hfill
  \begin{subfigure}[t]{0.4\linewidth}
\centering
\includegraphics[width=\linewidth]{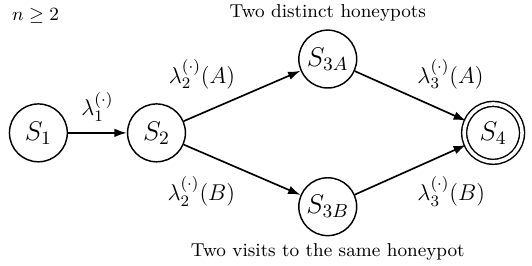}
\caption{$K=3$}
\label{fig:CTMCModelK3}
\end{subfigure}
  \caption{CTMC model for $K=2$ and $K=3$}
  \label{fig:service-level-param}
\end{figure}

We illustrate the state structure and transition dynamics for the case $K = 2$ in Fig.~\ref{fig:CTMCModelK2}. The CTMC consists of three states, $S_1$, $S_2$, and $S_3$, corresponding to the benign, suspicious, and malicious classifications, respectively. We begin by defining the transition rates that govern movement between these states. Let $n\ge 1$ be the integer number of honeypots deployed in the system. For each insider type $(\cdot) \in \{b,m\}$, we define two per-honeypot interaction rates. The parameter $\lambda^{(\cdot)}>0$ is the per-honeypot rate of a first interaction within the monitoring horizon. The parameter $\lambda^{\prime(\cdot)}\ge 0$ is the rate of a second interaction with a honeypot already visited once within that horizon. The superscripts $b$ and $m$ correspond to benign and malicious insiders, respectively, and capture differences in interaction intensity across insider types. For a given insider type, honeypots share the same per-visit rate parameters, and their interaction processes are assumed independent. Under these assumptions, the overall transition rate from state $S_1$ to state $S_2$ for a given insider type is
\begin{align}
    \lambda_1^{(\cdot)} = n \lambda^{(\cdot)},
\end{align}
since all $n$ honeypots are initially unvisited and each contributes independently to the occurrence of the first interaction event. After the first honeypot interaction has occurred and the insider enters state $S_2$, exactly one honeypot has been visited, while the remaining $n - 1$ honeypots are still unvisited. The overall transition rate from state $S_2$ to state $S_3$ is therefore given by
\begin{align}
    \lambda_2^{(\cdot)} = (n - 1)\lambda^{(\cdot)} + \lambda^{\prime(\cdot)},
\end{align}
where the first term captures interactions with any of the remaining unvisited honeypots, and the second term captures a repeated interaction with the previously visited honeypot. This formulation allows repeated interactions to occur at a rate that may differ from first-time interactions, thereby distinguishing first visits from repeated engagement while preserving the Markovian structure.

With the state ordering $(S_1, S_2, S_3)$, the CTMC generator matrices for benign and malicious insiders, denoted by $Q^{b}$ and $Q^{m}$, are upper triangular with an absorbing third state. For either insider type $(\cdot) \in \{b,m\}$, the generator matrix takes the form
\[
Q^{(\cdot)}
=
\begin{bmatrix}
-\lambda_1^{(\cdot)} & \lambda_1^{(\cdot)} & 0 \\
0 & -\lambda_2^{(\cdot)} & \lambda_2^{(\cdot)} \\
0 & 0 & 0
\end{bmatrix}
=
\begin{bmatrix}
-n \lambda^{(\cdot)}
& n \lambda^{(\cdot)}
& 0 \\[4pt]
0
& -\big((n-1)\lambda^{(\cdot)} + \lambda^{\prime (\cdot)}\big)
& (n-1)\lambda^{(\cdot)} + \lambda^{\prime (\cdot)} \\[4pt]
0 & 0 & 0
\end{bmatrix}.
\]

Let $\mathcal{P}^{(\cdot)}(T)=e^{Q^{(\cdot)}T}$ denote the transition probability matrix over the monitoring horizon $T$. If the CTMC begins in state $S_1$ at time 0, which is consistent with the assumption that all insiders are initially classified as benign, then the probability distribution at time $T$ is
\[
\mathbf{P}^{(\cdot)}(T)
=
[P^{(\cdot)}_1(T),\, P^{(\cdot)}_2(T),\, P^{(\cdot)}_3(T)]
=
\mathbf{e}_1^{\mathsf{T}} \mathcal{P}^{(\cdot)}(T),
\]
where $\mathbf{e}_1 = [1,0,0]^{\mathsf{T}}$. Because the chain has a series structure $S_1 \rightarrow S_2 \rightarrow S_3$ and the generator matrix is upper triangular, closed-form transient probabilities are available. For either insider type, when $\lambda^{(\cdot)}_1\ne\lambda^{(\cdot)}_2$, we have
\begin{align}
P^{(\cdot)}_1(T) &= e^{-n\,\lambda^{(\cdot)} T}, \nonumber \\[6pt]
P^{(\cdot)}_2(T) &=
\frac{n\,\lambda^{(\cdot)}}
{\big((n-1)\lambda^{(\cdot)} + \lambda^{\prime (\cdot)}\big) - n\lambda^{(\cdot)}}
\left( e^{-n\lambda^{(\cdot)}T}
- e^{-((n-1)\lambda^{(\cdot)} + \lambda^{\prime (\cdot)})T} \right), \nonumber \\[6pt]
P^{(\cdot)}_3(T) &= 1 - P^{(\cdot)}_1(T) - P^{(\cdot)}_2(T). \nonumber
\end{align}

When $\lambda^{(\cdot)}_1=\lambda^{(\cdot)}_2=r>0$, the limiting expressions are
\[
P^{(\cdot)}_1(T)=e^{-rT},\qquad
P^{(\cdot)}_2(T)=rT e^{-rT},\qquad
P^{(\cdot)}_3(T)=1-e^{-rT}(1+rT).
\]

Among these probabilities, $P^{(\cdot)}_3(T)$ represents the likelihood that the insider reaches the Malicious state within the detection horizon. For a malicious insider, $P^{ m}_3(T)$ quantifies the true-detection probability, whereas for a benign insider, $P^{b}_3(T)$ reflects the false-detection probability. These quantities provide a quantitative basis for evaluating the effectiveness of the detection framework in distinguishing malicious insiders from benign insiders.

The conditional expected detection time measures the time from the start of monitoring until the detection threshold is reached, among insiders detected within the horizon. Because the detection process is evaluated over a finite horizon of length $T$, the time to absorption is right-censored at $T$. As a result, the metric $\mathbb{E}[\tau^{(\cdot)}\mid\tau^{(\cdot)}\le T]$ conditions on detection by $T$ and differs from the unconditional mean censored time $\mathbb{E}[\min(\tau^{(\cdot)},T)]$.

For the two-stage model, the CTMC has a unidirectional structure $S_1 \rightarrow S_2 \rightarrow S_3$, where $S_3$ is the absorbing Malicious state. When $\lambda^{(\cdot)}_2>0$, the time to absorption $\tau^{(\cdot)}$ starting from $S_1$ is the sum of two independent exponentially distributed sojourn times with rates $\lambda^{(\cdot)}_1 = n\lambda^{(\cdot)}$ and $\lambda^{(\cdot)}_2 = (n-1)\lambda^{(\cdot)} + \lambda^{\prime (\cdot)}$. This yields a hypoexponential distribution for $\tau^{(\cdot)}$. For distinct positive rates, the expected detection time conditional on detection within the horizon is
\begin{equation}
\label{eq:K2_conditional_expected_detection_time}
\mathbb{E}\!\left[\tau^{(\cdot)} \mid \tau^{(\cdot)} \le T\right]
=
\frac{
\displaystyle
\frac{1}{\lambda^{(\cdot)}_2-\lambda^{(\cdot)}_1}
\left[
\frac{\lambda^{(\cdot)}_2}{\lambda^{(\cdot)}_1}
\left(1-e^{-\lambda^{(\cdot)}_1 T}\right)
-
\frac{\lambda^{(\cdot)}_1}{\lambda^{(\cdot)}_2}
\left(1-e^{-\lambda^{(\cdot)}_2 T}\right)
\right]
-
T\,\frac{
\lambda^{(\cdot)}_2 e^{-\lambda^{(\cdot)}_1 T}
-
\lambda^{(\cdot)}_1 e^{-\lambda^{(\cdot)}_2 T}
}{
\lambda^{(\cdot)}_2-\lambda^{(\cdot)}_1
}
}{
\displaystyle
1-
\frac{
\lambda^{(\cdot)}_2 e^{-\lambda^{(\cdot)}_1 T}
-
\lambda^{(\cdot)}_1 e^{-\lambda^{(\cdot)}_2 T}
}{
\lambda^{(\cdot)}_2-\lambda^{(\cdot)}_1
}
}.
\end{equation}

The derivation of \eqref{eq:K2_conditional_expected_detection_time} is provided in the Appendix. If $\lambda^{(\cdot)}_1=\lambda^{(\cdot)}_2=r>0$, the detection time has an Erlang distribution with two stages, and
\[
\mathbb{E}[\tau^{(\cdot)}\mid\tau^{(\cdot)}\le T]
=\frac{2}{r}\,
\frac{1-e^{-rT}\big(1+rT+(rT)^2/2\big)}
{1-e^{-rT}(1+rT)}.
\]
If $n=1$ and $\lambda^{\prime(\cdot)}=0$, the second interaction cannot occur, so $P^{(\cdot)}_3(T)=0$ and the conditional expected detection time is undefined. Whenever both stage rates are positive, the conditional mean converges to $1/\lambda^{(\cdot)}_1+1/\lambda^{(\cdot)}_2$ as $T\to\infty$.

\subsubsection{Three-stage Model}\label{Sec:Three-Step CTMC}

For $K=3$, detection requires three interactions, as illustrated in Fig.~\ref{fig:CTMCModelK3}. The detection levels follow the progression $(S_1,S_2,S_3,S_4)$, with $S_4$ the absorbing detection state. After two interactions, the next interaction rate depends on whether those visits involved distinct honeypots or the same honeypot. We therefore split $S_3$ into $S_{3A}$ for two distinct honeypots (Case A) and $S_{3B}$ for two visits to the same honeypot (Case B). For $n\ge 2$, the resulting state representation is $\{S_1,S_2,S_{3A},S_{3B},S_4\}$. It retains the visit history needed for the Markov property. For each insider type $(\cdot)\in\{b,m\}$, let $\lambda^{\prime\prime(\cdot)}\ge 0$ denote the per-honeypot rate of a third interaction with a honeypot visited twice within the monitoring horizon. The first transition rate matches the two-stage model:
\begin{equation}
\lambda_1^{(\cdot)} = n\,\lambda^{(\cdot)}.
\end{equation}
From $S_2$, the chain branches according to whether the second interaction involves a new honeypot (Case A) or the previously visited honeypot (Case B). In Case A, the transition rate is
\begin{equation}
\lambda_2^{(\cdot)}(A) = (n-1)\lambda^{(\cdot)}.
\end{equation}
In Case B, the transition rate is
\begin{equation}
\lambda_2^{(\cdot)}(B) = \lambda^{\prime(\cdot)}.
\end{equation}

From $S_{3A}$, the insider has interacted with two distinct honeypots, leaving $n-2$ unvisited honeypots available for a first interaction at rate $\lambda^{(\cdot)}$, while each of the two visited honeypots contributes a possible second interaction at rate $\lambda^{\prime(\cdot)}$. Therefore, the transition rate from $S_{3A}$ to $S_4$ is
\begin{equation}
\lambda_3^{(\cdot)}(A) = (n-2)\lambda^{(\cdot)} + 2\lambda^{\prime(\cdot)}.
\end{equation}
From $S_{3B}$, both prior interactions have occurred on the same honeypot, leaving $n-1$ unvisited honeypots available for first interactions at rate $\lambda^{(\cdot)}$, while the twice-visited honeypot contributes a possible third interaction at rate $\lambda^{\prime\prime(\cdot)}$. Hence, the transition rate from $S_{3B}$ to $S_4$ is
\begin{equation}
\lambda_3^{(\cdot)}(B) = (n-1)\lambda^{(\cdot)} + \lambda^{\prime\prime(\cdot)}.
\end{equation}

For $n\ge 2$, we write the CTMC generator matrix for each insider type $(\cdot)\in\{b,m\}$ as
\begin{gather}
Q^{(\cdot)}
=
\begin{bmatrix}
-\lambda_1^{(\cdot)} & \lambda_1^{(\cdot)} & 0 & 0 & 0 \\
0 & -\big(\lambda_2^{(\cdot)}(A)+\lambda_2^{(\cdot)}(B)\big) & \lambda_2^{(\cdot)}(A) & \lambda_2^{(\cdot)}(B) & 0 \\
0 & 0 & -\lambda_3^{(\cdot)}(A) & 0 & \lambda_3^{(\cdot)}(A) \\
0 & 0 & 0 & -\lambda_3^{(\cdot)}(B) & \lambda_3^{(\cdot)}(B) \\
0 & 0 & 0 & 0 & 0
\end{bmatrix}
\nonumber \\
=
\begin{bmatrix}
-n\lambda^{(\cdot)} & n\lambda^{(\cdot)} & 0 & 0 & 0\\
0 & -\big((n-1)\lambda^{(\cdot)} + \lambda^{\prime(\cdot)}\big)
  & (n-1)\lambda^{(\cdot)}  &  \lambda^{\prime(\cdot)} & 0\\
0 & 0 & -\big((n-2)\lambda^{(\cdot)} + 2\lambda^{\prime(\cdot)}\big) & 0
  & (n-2)\lambda^{(\cdot)} + 2\lambda^{\prime(\cdot)} \\
0 & 0 & 0 & -\big((n-1)\lambda^{(\cdot)} + \lambda^{\prime\prime(\cdot)}\big)
  & (n-1)\lambda^{(\cdot)} + \lambda^{\prime\prime(\cdot)} \\
0 & 0 & 0 & 0 & 0
\end{bmatrix}. \nonumber
\end{gather}

When $n=1$, Case A is impossible. The state $S_{3A}$ and its row and column must be omitted, leaving the chain $S_1\to S_2\to S_{3B}\to S_4$ with successive rates $\lambda^{(\cdot)}$, $\lambda^{\prime(\cdot)}$, and $\lambda^{\prime\prime(\cdot)}$. In particular, the expression for $\lambda^{(\cdot)}_3(A)$ must not be used for $n=1$, since it can be negative. The formulas below apply to the reduced chain by setting $P^{(\cdot)}_{3A}(T)=0$ and omitting all Case A terms.

For $n\ge 2$, given the transition probability matrix over the detection horizon $T$ as $\mathcal{P}^{(\cdot)}(T) = e^{Q^{(\cdot)}T}$, we have the transient probability vector at time $T$ for the three-stage model as
\[
\mathbf{P}^{(\cdot)}(T)
=
\big[P^{(\cdot)}_1(T),\,P^{(\cdot)}_2(T),\,P^{(\cdot)}_{3A}(T),\,P^{(\cdot)}_{3B}(T),\,P^{(\cdot)}_4(T)\big]
=
\mathbf{e}_1^{\mathsf{T}}\,\mathcal{P}^{(\cdot)}(T),
\]
where $\mathbf{e}_1=[1,0,0,0,0]^{\mathsf{T}}$. For $n\ge 2$, the upper-triangular generator yields the following transient probabilities when the rates appearing in each denominator are distinct. Coincident-rate cases are obtained by taking limits of the complete expressions or by evaluating $e^{Q^{(\cdot)}T}$:
\begin{align}
P^{(\cdot)}_1(T)
&=
e^{-\lambda^{(\cdot)}_1 T},
\nonumber\\[6pt]
P^{(\cdot)}_2(T)
&=
\frac{\lambda^{(\cdot)}_1}{\lambda^{(\cdot)}_2-\lambda^{(\cdot)}_1}
\left(
e^{-\lambda^{(\cdot)}_1 T}-e^{-\lambda^{(\cdot)}_2 T}
\right),
\nonumber\\[6pt]
P^{(\cdot)}_{3A}(T)
&=
\frac{\lambda^{(\cdot)}_2(A)\,\lambda^{(\cdot)}_1}{\lambda^{(\cdot)}_2-\lambda^{(\cdot)}_1}
\left[
\frac{e^{-\lambda^{(\cdot)}_3(A)T}-e^{-\lambda^{(\cdot)}_1T}}{\lambda^{(\cdot)}_1-\lambda^{(\cdot)}_3(A)}
-
\frac{e^{-\lambda^{(\cdot)}_3(A)T}-e^{-\lambda^{(\cdot)}_2T}}{\lambda^{(\cdot)}_2-\lambda^{(\cdot)}_3(A)}
\right],
\nonumber\\[6pt]
P^{(\cdot)}_{3B}(T)
&=
\frac{\lambda^{(\cdot)}_2(B)\,\lambda^{(\cdot)}_1}{\lambda^{(\cdot)}_2-\lambda^{(\cdot)}_1}
\left[
\frac{e^{-\lambda^{(\cdot)}_3(B)T}-e^{-\lambda^{(\cdot)}_1T}}{\lambda^{(\cdot)}_1-\lambda^{(\cdot)}_3(B)}
-
\frac{e^{-\lambda^{(\cdot)}_3(B)T}-e^{-\lambda^{(\cdot)}_2T}}{\lambda^{(\cdot)}_2-\lambda^{(\cdot)}_3(B)}
\right],
\nonumber\\[6pt]
P^{(\cdot)}_4(T)
&=
1-
P^{(\cdot)}_1(T)-P^{(\cdot)}_2(T)-P^{(\cdot)}_{3A}(T)-P^{(\cdot)}_{3B}(T),
\label{eq:K3_transient_probabilities}
\end{align}
where $\lambda^{(\cdot)}_2=\lambda^{(\cdot)}_2(A)+\lambda^{(\cdot)}_2(B)$. We present the details of deriving these equations in the Appendix section. In particular, the quantity $P^{(\cdot)}_4(T)$ gives the overall probability that an insider of type $(\cdot)$ reaches the Malicious state within the detection horizon. As before, $P^{ m}_4(T)$ represents the true-detection probability for a malicious insider, while $P^{b}_4(T)$ quantifies the false-detection probability for a benign insider under the three-stage model.

As in the two-stage model, we evaluate the expected detection time conditional on detection within the monitoring horizon. This quantity is defined only when $P^{(\cdot)}_4(T)>0$.

When $\lambda^{(\cdot)}_2>0$, the branch probabilities at the second interaction are
\[
\lambda^{(\cdot)}_2 = (n-1)\lambda^{(\cdot)}+\lambda^{\prime(\cdot)},
\qquad
p_A^{(\cdot)}=\frac{\lambda^{(\cdot)}_2(A)}{\lambda^{(\cdot)}_2}
=\frac{(n-1)\lambda^{(\cdot)}}{\lambda^{(\cdot)}_2},
\qquad
p_B^{(\cdot)}=\frac{\lambda^{(\cdot)}_2(B)}{\lambda^{(\cdot)}_2}
=\frac{\lambda^{\prime(\cdot)}}{\lambda^{(\cdot)}_2}.
\]
The branch choice is independent of the first two holding times. Conditional on a branch $j\in\{A,B\}$ with $p_j^{(\cdot)}>0$ and $\lambda^{(\cdot)}_3(j)>0$, the absorption time is a sum of independent exponential holding times:
\[
\tau^{(\cdot)}\mid j \;\;\overset{d}{=}\;\; X^{(\cdot)}_1 + X^{(\cdot)}_2 + X^{(\cdot)}_{3j},
\]
where $X^{(\cdot)}_1\sim \mathrm{Exp}(\lambda^{(\cdot)}_1)$, $X^{(\cdot)}_2\sim \mathrm{Exp}(\lambda^{(\cdot)}_2)$, and $X^{(\cdot)}_{3j}\sim \mathrm{Exp}(\lambda^{(\cdot)}_3(j))$, with independence. When every branch with positive probability has a positive final-stage rate, the survival function of $\tau^{(\cdot)}$ admits a closed form as a mixture of the two conditional survival functions:
\begin{equation}
\label{eq:K3_survival_mixture}
\mathbb{P}\!\left(\tau^{(\cdot)} > t\right)
=
p_A^{(\cdot)}\,\overline F_A^{(\cdot)}(t)\;+\;p_B^{(\cdot)}\,\overline F_B^{(\cdot)}(t),
\end{equation}
where, for a branch with positive probability and distinct positive rates $\lambda^{(\cdot)}_1,\lambda^{(\cdot)}_2,\lambda^{(\cdot)}_3(j)$,
\begin{align}
\overline F_j^{(\cdot)}(t)
&=
c_{1j}^{(\cdot)}\,e^{-\lambda^{(\cdot)}_1 t}
+
c_{2j}^{(\cdot)}\,e^{-\lambda^{(\cdot)}_2 t}
+
c_{3j}^{(\cdot)}\,e^{-\lambda^{(\cdot)}_3(j)\, t},
\label{eq:K3_survival_components}
\\ \nonumber
c_{1j}^{(\cdot)}
&=
\frac{\lambda^{(\cdot)}_2\,\lambda^{(\cdot)}_3(j)}
{\big(\lambda^{(\cdot)}_2-\lambda^{(\cdot)}_1\big)\big(\lambda^{(\cdot)}_3(j)-\lambda^{(\cdot)}_1\big)},
\quad
c_{2j}^{(\cdot)}
=
\frac{\lambda^{(\cdot)}_1\,\lambda^{(\cdot)}_3(j)}
{\big(\lambda^{(\cdot)}_1-\lambda^{(\cdot)}_2\big)\big(\lambda^{(\cdot)}_3(j)-\lambda^{(\cdot)}_2\big)},
\\ \nonumber
c_{3j}^{(\cdot)}
& =
\frac{\lambda^{(\cdot)}_1\,\lambda^{(\cdot)}_2}
{\big(\lambda^{(\cdot)}_1-\lambda^{(\cdot)}_3(j)\big)\big(\lambda^{(\cdot)}_2-\lambda^{(\cdot)}_3(j)\big)}.
\nonumber
\end{align}
In particular, $P^{(\cdot)}_4(T)=\mathbb{P}\!\left(\tau^{(\cdot)}\le T\right)=1-\mathbb{P}\!\left(\tau^{(\cdot)} > T\right)$, and $\mathbb{P}\!\left(\tau^{(\cdot)} > T\right)=P^{(\cdot)}_1(T)+P^{(\cdot)}_2(T)+P^{(\cdot)}_{3A}(T)+P^{(\cdot)}_{3B}(T)$, consistent with \eqref{eq:K3_transient_probabilities}.

Branches with zero probability are omitted from all mixture sums. If $\lambda^{(\cdot)}_2=0$, detection is impossible. If a branch with positive probability has $\lambda^{(\cdot)}_3(j)=0$, it never reaches the detection state, so its survival function is $\overline F_j^{(\cdot)}(t)=1$ and its contribution to integrated survival is $p_j^{(\cdot)}T$. These cases can also be evaluated directly using the appropriate generator. For coincident positive rates, limits must be taken for the complete survival expressions, not for individual coefficients that may diverge.

To obtain $\mathbb{E}\!\left[\tau^{(\cdot)} \mid \tau^{(\cdot)} \le T\right]$ in closed form, we use the identity
\[
\mathbb{E}\!\left[\tau^{(\cdot)}\,\mathbf{1}\{\tau^{(\cdot)}\le T\}\right]
=
\int_0^T \mathbb{P}\!\left(\tau^{(\cdot)} > t\right)\,dt
\;-\;
T\,\mathbb{P}\!\left(\tau^{(\cdot)} > T\right),
\]
which follows from $\mathbb{E}[\min(\tau^{(\cdot)},T)]=\int_0^T \mathbb{P}(\tau^{(\cdot)} > t)\,dt$. For the distinct positive-rate case, define
\begin{equation}
\label{eq:K3_integrated_survival}
\begin{aligned}
I^{(\cdot)}(T)
&:=
\int_0^T \mathbb{P}\!\left(\tau^{(\cdot)} > t\right)\,dt\\
&=
\sum_{j\in\{A,B\}} p_j^{(\cdot)}
\left[
\frac{c_{1j}^{(\cdot)}}{\lambda^{(\cdot)}_1}\Big(1-e^{-\lambda^{(\cdot)}_1 T}\Big)
+
\frac{c_{2j}^{(\cdot)}}{\lambda^{(\cdot)}_2}\Big(1-e^{-\lambda^{(\cdot)}_2 T}\Big)
+
\frac{c_{3j}^{(\cdot)}}{\lambda^{(\cdot)}_3(j)}\Big(1-e^{-\lambda^{(\cdot)}_3(j) T}\Big)
\right].
\end{aligned}
\end{equation}
For $P^{(\cdot)}_4(T)>0$, the expected detection time conditional on detection within the horizon is
\begin{equation}
\label{eq:K3_conditional_expected_detection_time}
\mathbb{E}\!\left[\tau^{(\cdot)} \mid \tau^{(\cdot)} \le T\right]
=
\frac{
I^{(\cdot)}(T)
-
T\,\mathbb{P}\!\left(\tau^{(\cdot)} > T\right)
}{
1-\mathbb{P}\!\left(\tau^{(\cdot)} > T\right)
}
=
\frac{
I^{(\cdot)}(T)
-
T\big(1-P^{(\cdot)}_4(T)\big)
}{
P^{(\cdot)}_4(T)
}.
\end{equation}
The derivation of \eqref{eq:K3_survival_mixture}--\eqref{eq:K3_conditional_expected_detection_time} is provided in the Appendix. The conditional-mean identity remains valid in the boundary cases when survival probabilities and their integrals are evaluated as described above and detection has positive probability. If $\lambda^{(\cdot)}_2>0$ and every branch with positive probability has a positive final-stage rate, detection occurs eventually with probability one. Under these conditions, as $T\to\infty$, \eqref{eq:K3_conditional_expected_detection_time} converges to
\[
\mathbb{E}[\tau^{(\cdot)}]
=
\frac{1}{\lambda^{(\cdot)}_1}
+
\frac{1}{\lambda^{(\cdot)}_2}
+
\sum_{\substack{j\in\{A,B\}\\p_j^{(\cdot)}>0}}
\frac{p_j^{(\cdot)}}{\lambda^{(\cdot)}_3(j)}.
\]
The first two holding-time distributions are common to both branches, while the final holding time depends on the visit history.

\subsection{Optimizing the Detection Process} \label{Sec:Optimization}

In the final step of the framework, the system operator selects the number of honeypots $n$ and the monitoring horizon $T$ to minimize the expected costs of malicious exposure and false detections. This use of detection measures to guide resource decisions is related to strategic inspection models such as that of \citet{dahan2022network}, who minimize the number of detectors needed to achieve a target expected detection rate against worst-case attacks. Here, the CTMC probabilities and detection-time measures from Section~\ref{Sec:CTMC} quantify the consequences of each configuration for both insider types. The objective accounts for damage incurred before malicious detection or after a missed detection, together with the investigation effort and operational disruption caused by falsely flagging benign insiders.

Let $Z\ge1$ be the number of employees with access to sensitive system information. We evaluate cost conditional on exactly one of these employees being malicious during the monitoring horizon, with the remaining $Z-1$ employees benign. The malicious insider's identity is unknown, and the same detection rule is applied to all employees. This single-insider scenario provides a tractable conditional assessment without requiring a probability model for the number of malicious employees or assumptions about coordination among them. Correct detection and threat removal are assumed to stop further damage from the detected insider. The false-detection cost evaluates each benign employee over the full monitoring horizon $T$.

Let $c^a\ge0$ denote the damage cost per unit time while the malicious insider remains undetected and retains access. We fix an exposure cutoff $\bar T>0$, independently of the decision $T$, and require $T\le\bar T$. The objective evaluates a single monitoring horizon and credits only malicious detections occurring by $T$. If detection occurs by $T$, damage accumulates until the detection time. Otherwise, the cost model assigns exposure from time $0$ through $\bar T$, representing containment outside the modeled detection process. Assigning the full duration $\bar T$ to every missed detection is a conservative convention when such containment could occur earlier. Let $c^f\ge0$ denote the expected cost of falsely detecting one benign insider, including investigation effort, operational disruption, and productivity loss. We define the resulting expected cost as
\begin{equation}
C(n,T)
=c^a\left[
\int_0^T t\,\frac{\mathrm{d}P^m_{K+1}(t)}{\mathrm{d} t} \mathrm{d} t
+\big(1-P^m_{K+1}(T)\big)\bar T
\right]
+(Z-1)c^fP^b_{K+1}(T).
\label{eq:defender_objective_time}
\end{equation}
where the first term captures the expected attack damage accumulated over the malicious insider's exposure period, and the second term captures the expected cost of false detections among benign insiders. When $P^m_{K+1}(T)>0$, the integral equals
$P^m_{K+1}(T)\mathbb{E}[\tau^m\mid\tau^m\le T]$.
Writing this term as an integral also covers configurations with zero detection probability, for which it is zero and the conditional mean is undefined. The first cost component is the expected damage under the exposure convention above. The second follows by summing the expected false-detection costs over the $Z-1$ benign employees. This sum does not require independence across employees.

Let $\mathcal N$ be a nonempty finite set of feasible positive integer honeypot counts, and let $\mathcal T$ be a nonempty compact subset of $(0,\bar T]$, such as a finite grid or a closed interval of permitted monitoring horizons. For a fixed detection threshold $K$, the optimization problem is
\begin{equation}
\min_{n\in\mathcal N,\,T\in\mathcal T} C(n,T).
\label{eq:opt_problem}
\end{equation}

Proposition~\ref{prop:marginal_monitoring_value} characterizes the marginal value of monitoring. Let $f_m(t)=\mathrm{d}P^m_{K+1}(t)/\mathrm{d}t$ and $f_b(t)=\mathrm{d}P^b_{K+1}(t)/\mathrm{d}t$. A small extension of the horizon admits additional malicious detections that replace the assigned exposure $\bar T$ with detection times near $T$, yielding a marginal saving of $c^a(\bar T-T)f_m(T)$. The competing term, $(Z-1)c^f f_b(T)$, captures the marginal cost of additional false detections. Their balance depends on deployment size and interaction rates, linking the monitoring decision to repeated honeypot engagement. The equality condition alone establishes neither optimality nor uniqueness. For a finite grid of permitted horizons, the decision follows from comparing their expected costs. The proof is provided in Appendix~\ref{app:cost_derivation}.

\begin{proposition}
\label{prop:marginal_monitoring_value}
Fix $n$, $K$, the interaction-rate parameters, and the exposure cutoff $\bar T$. For $0<T<\bar T$,
\begin{equation}
\frac{\partial C(n,T)}{\partial T}
=-c^a(\bar T-T)f_m(T)+(Z-1)c^f f_b(T).
\label{eq:cost_horizon_derivative}
\end{equation}
Consequently, the marginal cost of extending monitoring is negative when
$c^a(\bar T-T)f_m(T)>(Z-1)c^f f_b(T)$.
Any interior minimizer over a continuous feasible interval must satisfy equality between these two terms.
\end{proposition}

\section{Structural Properties of Detection Performance}\label{sec:TheoreticalProperties}

This section establishes structural properties of the CTMC model for the two-stage ($K=2$) and three-stage ($K=3$) configurations. We identify combinations of honeypot counts and repeated-interaction rates that preclude false detection and establish conditions under which malicious insiders are detected at least as readily as benign insiders. Specifically, when malicious interaction rates are no smaller at each per-honeypot visit level, malicious detection times are stochastically no larger and detection probabilities are at least as high for any common honeypot count and monitoring horizon. For $K=2$, the expected detection time conditional on detection within the horizon is also no greater for malicious insiders, whenever both conditional means are defined. We further show that increasing the detection threshold from $K=2$ to $K=3$, with other parameters fixed, weakly reduces both true- and false-detection probabilities. These results clarify the trade-off between detecting malicious activity and avoiding false detections and provide a basis for interpreting the subsequent numerical results. All proofs are provided in Appendix~\ref{app:analytical_results}.

\subsection{Structural Analysis of Detection Performance}\label{sec:theory_K2}

Consider the case $K=2$, corresponding to the two-stage model described in Section~\ref{Sec:Two-Step CTMC}. For benign insiders, repeated engagement with the same honeypot may be limited, leading to smaller values of $\lambda^{\prime b}$ relative to the rate of first-time interactions. Malicious insiders, in contrast, may deliberately revisit or pursue deeper interaction with honeypots they have already triggered, resulting in higher interaction intensities across both stages of the detection process. Motivated by these behavioral differences, we evaluate the following cross-type rate dominance conditions:
\begin{equation}\label{eq:rate_assumptions_K2_prop}
\lambda^{m}\ge \lambda^{b},\qquad
\lambda^{\prime m}\ge \lambda^{\prime b}.
\end{equation}
Further, these conditions ensure that the aggregate transition rate out of each transient state in the two-stage model is no smaller for malicious insiders than for benign insiders.

We first consider the boundary case $\lambda^{\prime b}=0$, in which a benign insider interacts with each honeypot at most once. The following result identifies how this restriction affects false detection when only one honeypot is deployed.

\begin{lemma}\label{lem:n1_no_false_detection}
For $K=2$, $n=1$, and $\lambda^{\prime b}=0$, the false-detection probability satisfies $P^b_3(T)=0$ for every $T\ge0$.
\end{lemma}

Under this configuration, benign insiders cannot accumulate the two interactions required for detection. Malicious detection remains possible if $\lambda^{\prime m}>0$, but is also impossible if $\lambda^{\prime m}=0$. Thus, eliminating false detections does not by itself establish an effective detection configuration. Its value also depends on the rate at which malicious insiders supply the required repeated interaction.

The next proposition compares the detection-time distributions under the rate conditions in~\eqref{eq:rate_assumptions_K2_prop}. This property implies that the detection mechanism escalates malicious behavior no more slowly than it escalates benign behavior, thereby supporting both detection accuracy and responsiveness.

\begin{proposition}\label{Prop:DetectionTime}
Consider the two-stage model with a common honeypot count $n\ge1$ and the rate conditions in \eqref{eq:rate_assumptions_K2_prop}. Then
\begin{equation}
\tau^m\le_{\mathrm{st}}\tau^b,
\qquad
P^m_3(T)\ge P^b_3(T)\quad\text{for every }T\ge0.
\label{eq:K2_stochastic_order}
\end{equation}
In particular, the two detection times admit a coupling under which $\tau^m\le\tau^b$ almost surely. If $T>0$ and $(n-1)\lambda^b+\lambda^{\prime b}>0$, then both conditional means are defined and
\begin{equation}
\mathbb{E}[\tau^m\mid\tau^m\le T]
\le
\mathbb{E}[\tau^b\mid\tau^b\le T].
\label{eq:K2_conditional_order}
\end{equation}
\end{proposition}

Appendix~\ref{app:proof_K2_order} proves the probability comparison by coupling the holding times and establishes the conditional comparison separately. The latter requires an additional argument because conditioning each process on detection by $T$ selects different sets of realizations. Since $\lambda^b>0$, its positivity condition holds for every $n\ge2$, and also for $n=1$ when $\lambda^{\prime b}>0$.

\subsection{Analytical Properties of the Detection Models}\label{sec:theory_K3}

We next examine how increasing the detection threshold to $K=3$ affects the relative detection behavior of benign and malicious insiders. Although the three-stage model introduces history-dependent transition dynamics through the branching at the second interaction, it is still possible to establish a rigorous ordering of detection outcomes under mild rate conditions.

Case A records visits to two distinct honeypots, whereas Case B records two visits to the same honeypot. To extend the cross-type comparison to this setting, we assume
\begin{equation}\label{eq:rate_assumptions_K3_prop}
\lambda^m\ge\lambda^b,\qquad
\lambda^{\prime m}\ge\lambda^{\prime b},\qquad
\lambda^{\prime\prime m}\ge\lambda^{\prime\prime b}.
\end{equation}
These conditions require malicious insiders to interact at least as frequently at each per-honeypot visit level, as formalized in Proposition~\ref{Prop:DetectionTimeK3}. The branch probabilities also depend on the rates, so the proof compares the underlying per-honeypot interaction counts while allowing the two insider types to follow different histories.

\begin{proposition}\label{Prop:DetectionTimeK3}
Consider the three-stage model with a common honeypot count $n\ge1$ and the rate conditions in \eqref{eq:rate_assumptions_K3_prop}. For $n=1$, use the reduced chain specified in Section~\ref{Sec:Three-Step CTMC}. Then
\begin{equation}
\tau^m\le_{\mathrm{st}}\tau^b,
\qquad
P^m_4(T)\ge P^b_4(T)\quad\text{for every }T\ge0.
\label{eq:P4_ordering_K3_prop}
\end{equation}
The two detection times admit a coupling under which $\tau^m\le\tau^b$ almost surely.
\end{proposition}

The coupling in Appendix~\ref{app:proof_K3_order} preserves the marginal interaction rates of each insider type and orders their cumulative interaction counts. This proves the detection-probability comparison without requiring the two processes to take the same branch.

False detection is possible only if a benign insider can complete three interactions. Because $\lambda^b>0$, the following conditions are necessary and sufficient for $P^b_4(T)>0$ at any $T>0$:
\begin{itemize}
\item If $\lambda^{\prime b}=0$, every interaction must involve a distinct honeypot, so $n\ge3$ is required.
\item If $\lambda^{\prime b}>0$ and $\lambda^{\prime\prime b}=0$, a honeypot can be visited twice, and $n\ge2$ is required.
\item If $\lambda^{\prime b}>0$ and $\lambda^{\prime\prime b}>0$, a single honeypot suffices, so every $n\ge1$ permits false detection.
\end{itemize}
The path arguments are given in Appendix~\ref{app:proof_K3_order}. A configuration with zero false-detection probability can still detect malicious insiders if their repeated-interaction rates permit absorption. Such configurations therefore remain eligible for optimization. The associated benign conditional detection time is undefined, rather than zero.

Increasing the threshold from $K=2$ to $K=3$ requires an additional interaction before detection. To isolate this effect, fix $n$ and the interaction rates and apply both thresholds to the same underlying interaction process. The third interaction cannot precede the second, so
\[
\tau^{(\cdot)}_{K=2}\le\tau^{(\cdot)}_{K=3},
\qquad
\{\tau^{(\cdot)}_{K=3}\le T\}\subseteq
\{\tau^{(\cdot)}_{K=2}\le T\}.
\]
Consequently,
\begin{equation}
P^b_4(T)\le P^b_3(T),\qquad
P^m_4(T)\le P^m_3(T),
\label{eq:threshold_probability_order}
\end{equation}
where the subscript $3$ denotes the absorbing state of the two-stage model and the subscript $4$ that of the three-stage model. Thus, requiring more evidence weakly reduces both true- and false-detection probabilities within a fixed horizon. This comparison explains the trade-off between timely malicious detection and protection against false detections. It does not, by itself, order the expected times conditional on detection under the two different thresholds.

\section{Experimental Results}\label{sec:experiments}

This section evaluates the proposed CTMC-based insider detection framework using the analytical measures derived in Section~\ref{Sec:CTMC}, the cost function in Section~\ref{Sec:Optimization}, and the structural properties in Section~\ref{sec:TheoreticalProperties}. We examine how transient state probabilities, true- and false-detection probabilities, expected detection times conditional on detection within the horizon, and expected total cost vary with the honeypot count $n$, monitoring horizon $T$, and detection threshold $K$. %

\subsection{Experiment Setup}\label{Sec:ExprimentSetup}

We adopt the cross-type rate-dominance assumptions in \eqref{eq:rate_assumptions_K2_prop} for $K=2$ and \eqref{eq:rate_assumptions_K3_prop} for $K=3$. These conditions specify scenarios in which malicious interaction rates are no smaller than benign interaction rates at each per-honeypot visit level. The parameter values below define illustrative scenarios for examining how interaction frequency and repeated engagement affect detection performance. %

For benign insiders, we select the first-interaction rate $\lambda^{b}$ from $\{0.1,0.2,0.3\}$ and set $\lambda^{\prime b}=0.05$, representing a small but positive rate of second interaction with the same honeypot. For $K=3$, we set $\lambda^{\prime\prime b}=0$, so a benign insider interacts with each honeypot at most twice. For malicious insiders, we consider first-interaction rates $\lambda^{m}\in\{0.5,1,2\}$, set $\lambda^{\prime m}=2$, and use third-interaction rates $\lambda^{\prime\prime m}\in\{0.5,4\}$ for $K=3$. The latter values lie below and above $\lambda^{\prime m}$, respectively, allowing us to compare slower and faster third interactions with the same honeypot. All selected rates satisfy the stated cross-type inequalities.

We consider honeypot counts $n\in\{1,\ldots,10\}$ and monitoring horizons $T\in\{2,4,6\}$. All times use a common time unit, and interaction rates are expressed per unit time. Under $K=3$, the configuration $n=1$ has zero false-detection probability because $\lambda^{\prime\prime b}=0$, as established in Section~\ref{sec:theory_K3}. Conditional expected detection times are evaluated only when the corresponding detection probability is positive. The cost parameters used in the optimization experiments are specified in Section~\ref{sec:theory_optimization}.

\subsection{Monotonicity and Sensitivity of Transient State Probabilities}
\label{sec:Experiment_theory_monotonicity}

\begin{figure*}[tbp]
\centering
\captionsetup[sub]{font=small}
\includegraphics[width=0.4\linewidth]{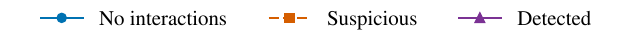}\par\smallskip
\begin{subfigure}[t]{0.27\linewidth}
\centering
\includegraphics[width=\linewidth]{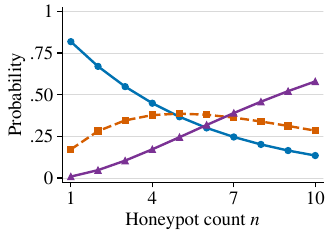}
\caption{$T=2$, $\lambda^b=0.1$}
\label{fig:k2-benign-t2-r1}
\end{subfigure}\hfill
\begin{subfigure}[t]{0.27\linewidth}
\centering
\includegraphics[width=\linewidth]{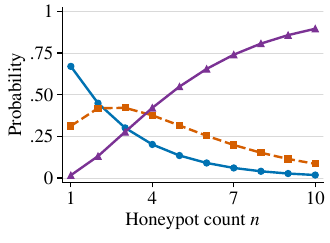}
\caption{$T=2$, $\lambda^b=0.2$}
\label{fig:k2-benign-t2-r2}
\end{subfigure}\hfill
\begin{subfigure}[t]{0.27\linewidth}
\centering
\includegraphics[width=\linewidth]{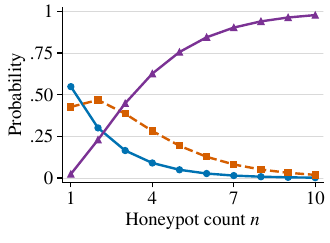}
\caption{$T=2$, $\lambda^b=0.3$}
\label{fig:k2-benign-t2-r3}
\end{subfigure}\par\medskip
\begin{subfigure}[t]{0.27\linewidth}
\centering
\includegraphics[width=\linewidth]{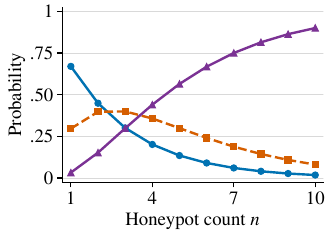}
\caption{$T=4$, $\lambda^b=0.1$}
\label{fig:k2-benign-t4-r1}
\end{subfigure}\hfill
\begin{subfigure}[t]{0.27\linewidth}
\centering
\includegraphics[width=\linewidth]{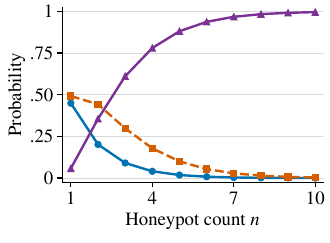}
\caption{$T=4$, $\lambda^b=0.2$}
\label{fig:k2-benign-t4-r2}
\end{subfigure}\hfill
\begin{subfigure}[t]{0.27\linewidth}
\centering
\includegraphics[width=\linewidth]{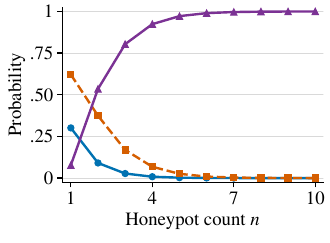}
\caption{$T=4$, $\lambda^b=0.3$}
\label{fig:k2-benign-t4-r3}
\end{subfigure}\par\medskip
\begin{subfigure}[t]{0.27\linewidth}
\centering
\includegraphics[width=\linewidth]{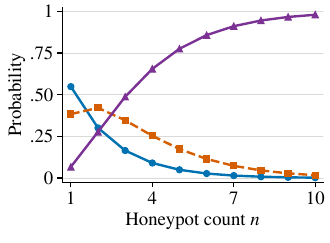}
\caption{$T=6$, $\lambda^b=0.1$}
\label{fig:k2-benign-t6-r1}
\end{subfigure}\hfill
\begin{subfigure}[t]{0.27\linewidth}
\centering
\includegraphics[width=\linewidth]{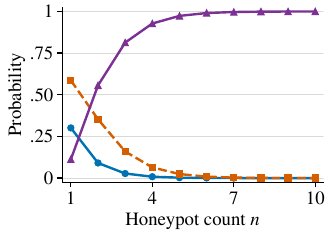}
\caption{$T=6$, $\lambda^b=0.2$}
\label{fig:k2-benign-t6-r2}
\end{subfigure}\hfill
\begin{subfigure}[t]{0.27\linewidth}
\centering
\includegraphics[width=\linewidth]{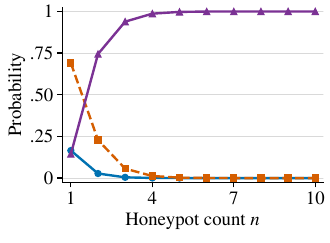}
\caption{$T=6$, $\lambda^b=0.3$}
\label{fig:k2-benign-t6-r3}
\end{subfigure}\par\medskip
\caption{Transient state probabilities for benign insiders with $K=2$ and $\lambda^{\prime b}=0.05$.}
\label{fig: Prob_MC_benign_K2}
\end{figure*}

\begin{figure*}[tbp]
\centering
\captionsetup[sub]{font=small}
\includegraphics[width=0.4\linewidth]{Figures/transient_probabilities/legend.pdf}\par\smallskip
\begin{subfigure}[t]{0.27\linewidth}
\centering
\includegraphics[width=\linewidth]{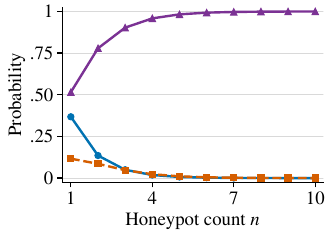}
\caption{$T=2$, $\lambda^m=0.5$}
\label{fig:k2-malicious-t2-r1}
\end{subfigure}\hfill
\begin{subfigure}[t]{0.27\linewidth}
\centering
\includegraphics[width=\linewidth]{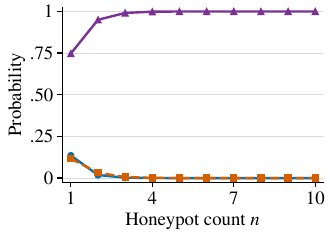}
\caption{$T=2$, $\lambda^m=1$}
\label{fig:k2-malicious-t2-r2}
\end{subfigure}\hfill
\begin{subfigure}[t]{0.27\linewidth}
\centering
\includegraphics[width=\linewidth]{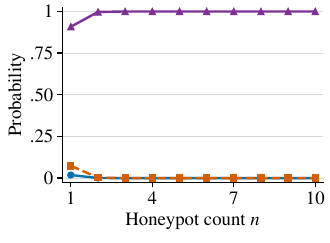}
\caption{$T=2$, $\lambda^m=2$}
\label{fig:k2-malicious-t2-r3}
\end{subfigure}\par\medskip
\begin{subfigure}[t]{0.27\linewidth}
\centering
\includegraphics[width=\linewidth]{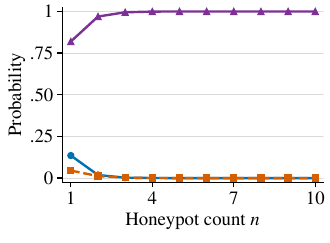}
\caption{$T=4$, $\lambda^m=0.5$}
\label{fig:k2-malicious-t4-r1}
\end{subfigure}\hfill
\begin{subfigure}[t]{0.27\linewidth}
\centering
\includegraphics[width=\linewidth]{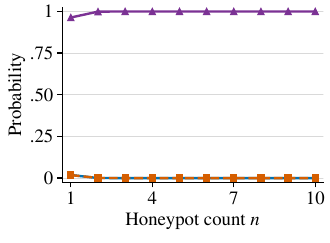}
\caption{$T=4$, $\lambda^m=1$}
\label{fig:k2-malicious-t4-r2}
\end{subfigure}\hfill
\begin{subfigure}[t]{0.27\linewidth}
\centering
\includegraphics[width=\linewidth]{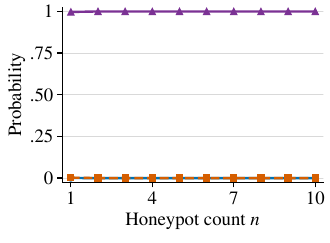}
\caption{$T=4$, $\lambda^m=2$}
\label{fig:k2-malicious-t4-r3}
\end{subfigure}\par\medskip
\begin{subfigure}[t]{0.27\linewidth}
\centering
\includegraphics[width=\linewidth]{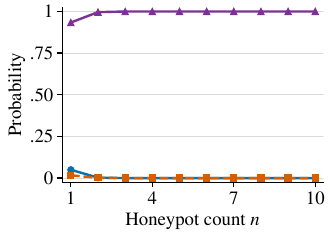}
\caption{$T=6$, $\lambda^m=0.5$}
\label{fig:k2-malicious-t6-r1}
\end{subfigure}\hfill
\begin{subfigure}[t]{0.27\linewidth}
\centering
\includegraphics[width=\linewidth]{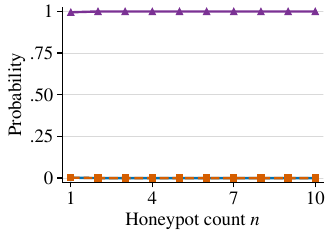}
\caption{$T=6$, $\lambda^m=1$}
\label{fig:k2-malicious-t6-r2}
\end{subfigure}\hfill
\begin{subfigure}[t]{0.27\linewidth}
\centering
\includegraphics[width=\linewidth]{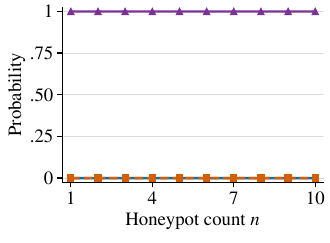}
\caption{$T=6$, $\lambda^m=2$}
\label{fig:k2-malicious-t6-r3}
\end{subfigure}\par\medskip
\caption{Transient state probabilities for malicious insiders with $K=2$ and $\lambda^{\prime m}=2$.}
\label{fig: Prob_MC_malicious_K2}
\end{figure*}

\begin{figure*}[tbp]
\centering
\captionsetup[sub]{font=small}
\includegraphics[width=0.4\linewidth]{Figures/transient_probabilities/legend.pdf}\par\smallskip
\begin{subfigure}[t]{0.27\linewidth}
\centering
\includegraphics[width=\linewidth]{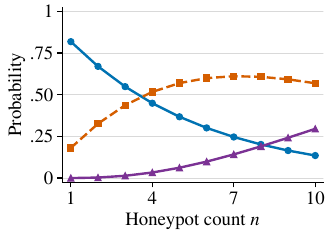}
\caption{$T=2$, $\lambda^b=0.1$}
\label{fig:k3-benign-t2-r1}
\end{subfigure}\hfill
\begin{subfigure}[t]{0.27\linewidth}
\centering
\includegraphics[width=\linewidth]{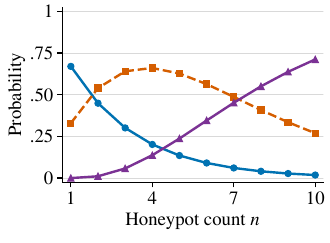}
\caption{$T=2$, $\lambda^b=0.2$}
\label{fig:k3-benign-t2-r2}
\end{subfigure}\hfill
\begin{subfigure}[t]{0.27\linewidth}
\centering
\includegraphics[width=\linewidth]{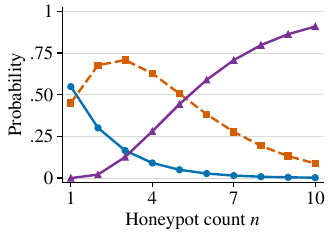}
\caption{$T=2$, $\lambda^b=0.3$}
\label{fig:k3-benign-t2-r3}
\end{subfigure}\par\medskip
\begin{subfigure}[t]{0.27\linewidth}
\centering
\includegraphics[width=\linewidth]{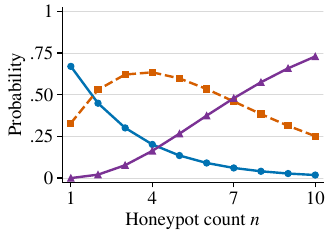}
\caption{$T=4$, $\lambda^b=0.1$}
\label{fig:k3-benign-t4-r1}
\end{subfigure}\hfill
\begin{subfigure}[t]{0.27\linewidth}
\centering
\includegraphics[width=\linewidth]{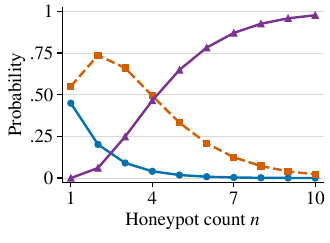}
\caption{$T=4$, $\lambda^b=0.2$}
\label{fig:k3-benign-t4-r2}
\end{subfigure}\hfill
\begin{subfigure}[t]{0.27\linewidth}
\centering
\includegraphics[width=\linewidth]{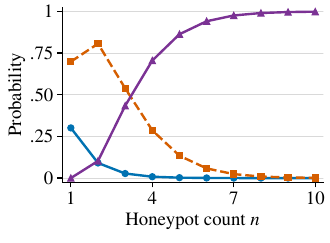}
\caption{$T=4$, $\lambda^b=0.3$}
\label{fig:k3-benign-t4-r3}
\end{subfigure}\par\medskip
\begin{subfigure}[t]{0.27\linewidth}
\centering
\includegraphics[width=\linewidth]{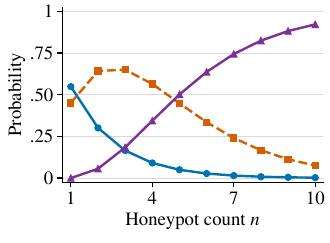}
\caption{$T=6$, $\lambda^b=0.1$}
\label{fig:k3-benign-t6-r1}
\end{subfigure}\hfill
\begin{subfigure}[t]{0.27\linewidth}
\centering
\includegraphics[width=\linewidth]{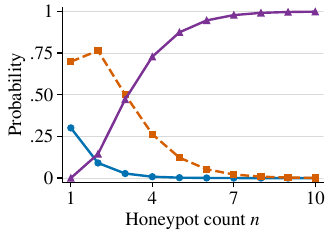}
\caption{$T=6$, $\lambda^b=0.2$}
\label{fig:k3-benign-t6-r2}
\end{subfigure}\hfill
\begin{subfigure}[t]{0.27\linewidth}
\centering
\includegraphics[width=\linewidth]{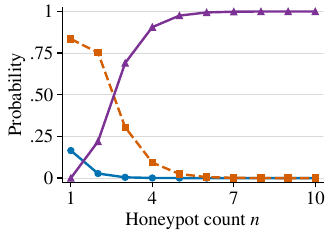}
\caption{$T=6$, $\lambda^b=0.3$}
\label{fig:k3-benign-t6-r3}
\end{subfigure}\par\medskip
\caption{Transient state probabilities for benign insiders with $K=3$, $\lambda^{\prime b}=0.05$, and $\lambda^{\prime\prime b}=0$.}
\label{fig: Prob_MC_benign_K3}
\end{figure*}

\begin{figure*}[tbp]
\centering
\captionsetup[sub]{font=small}
\includegraphics[width=0.4\linewidth]{Figures/transient_probabilities/legend.pdf}\par\smallskip
\begin{subfigure}[t]{0.27\linewidth}
\centering
\includegraphics[width=\linewidth]{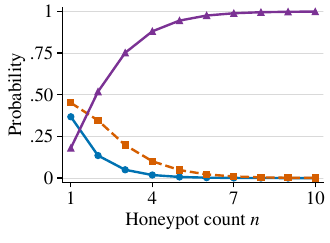}
\caption{$T=2$, $\lambda^m=0.5$}
\label{fig:k3-malicious-05-t2-r1}
\end{subfigure}\hfill
\begin{subfigure}[t]{0.27\linewidth}
\centering
\includegraphics[width=\linewidth]{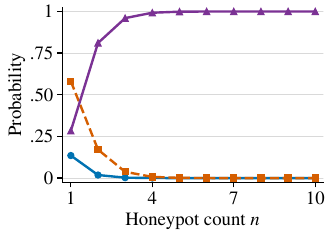}
\caption{$T=2$, $\lambda^m=1$}
\label{fig:k3-malicious-05-t2-r2}
\end{subfigure}\hfill
\begin{subfigure}[t]{0.27\linewidth}
\centering
\includegraphics[width=\linewidth]{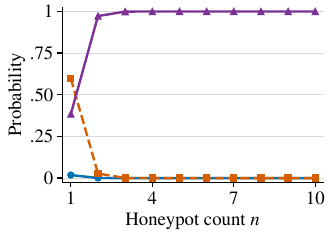}
\caption{$T=2$, $\lambda^m=2$}
\label{fig:k3-malicious-05-t2-r3}
\end{subfigure}\par\medskip
\begin{subfigure}[t]{0.27\linewidth}
\centering
\includegraphics[width=\linewidth]{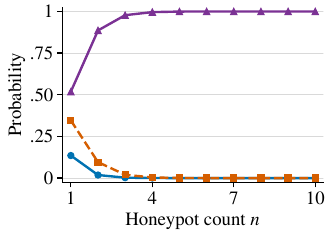}
\caption{$T=4$, $\lambda^m=0.5$}
\label{fig:k3-malicious-05-t4-r1}
\end{subfigure}\hfill
\begin{subfigure}[t]{0.27\linewidth}
\centering
\includegraphics[width=\linewidth]{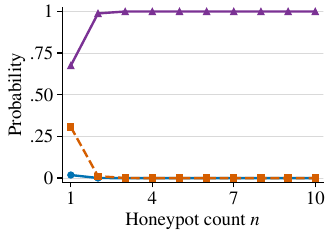}
\caption{$T=4$, $\lambda^m=1$}
\label{fig:k3-malicious-05-t4-r2}
\end{subfigure}\hfill
\begin{subfigure}[t]{0.27\linewidth}
\centering
\includegraphics[width=\linewidth]{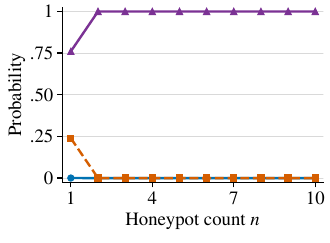}
\caption{$T=4$, $\lambda^m=2$}
\label{fig:k3-malicious-05-t4-r3}
\end{subfigure}\par\medskip
\begin{subfigure}[t]{0.27\linewidth}
\centering
\includegraphics[width=\linewidth]{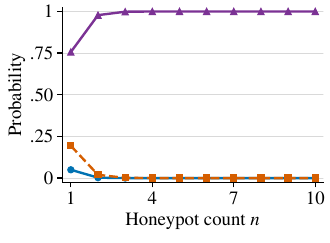}
\caption{$T=6$, $\lambda^m=0.5$}
\label{fig:k3-malicious-05-t6-r1}
\end{subfigure}\hfill
\begin{subfigure}[t]{0.27\linewidth}
\centering
\includegraphics[width=\linewidth]{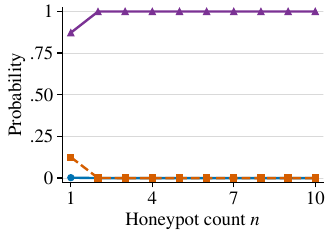}
\caption{$T=6$, $\lambda^m=1$}
\label{fig:k3-malicious-05-t6-r2}
\end{subfigure}\hfill
\begin{subfigure}[t]{0.27\linewidth}
\centering
\includegraphics[width=\linewidth]{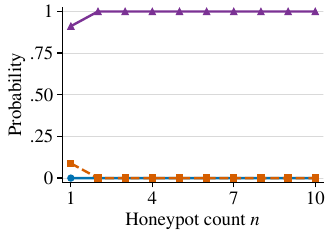}
\caption{$T=6$, $\lambda^m=2$}
\label{fig:k3-malicious-05-t6-r3}
\end{subfigure}\par\medskip
\caption{Transient state probabilities for malicious insiders with $K=3$, $\lambda^{\prime m}=2$, and $\lambda^{\prime\prime m}=0.5$.}
\label{fig: Prob_MC_malicious_K3_05}
\end{figure*}

\begin{figure*}[tbp]
\centering
\captionsetup[sub]{font=small}
\includegraphics[width=0.4\linewidth]{Figures/transient_probabilities/legend.pdf}\par\smallskip
\begin{subfigure}[t]{0.27\linewidth}
\centering
\includegraphics[width=\linewidth]{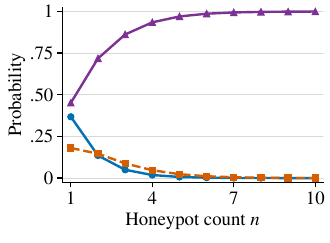}
\caption{$T=2$, $\lambda^m=0.5$}
\label{fig:k3-malicious-4-t2-r1}
\end{subfigure}\hfill
\begin{subfigure}[t]{0.27\linewidth}
\centering
\includegraphics[width=\linewidth]{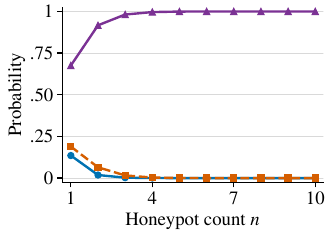}
\caption{$T=2$, $\lambda^m=1$}
\label{fig:k3-malicious-4-t2-r2}
\end{subfigure}\hfill
\begin{subfigure}[t]{0.27\linewidth}
\centering
\includegraphics[width=\linewidth]{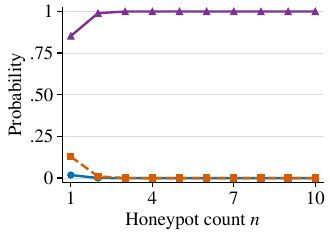}
\caption{$T=2$, $\lambda^m=2$}
\label{fig:k3-malicious-4-t2-r3}
\end{subfigure}\par\medskip
\begin{subfigure}[t]{0.27\linewidth}
\centering
\includegraphics[width=\linewidth]{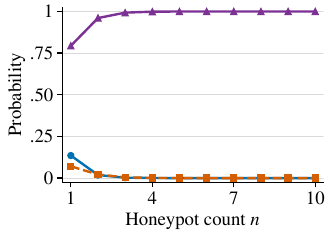}
\caption{$T=4$, $\lambda^m=0.5$}
\label{fig:k3-malicious-4-t4-r1}
\end{subfigure}\hfill
\begin{subfigure}[t]{0.27\linewidth}
\centering
\includegraphics[width=\linewidth]{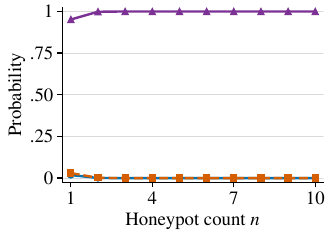}
\caption{$T=4$, $\lambda^m=1$}
\label{fig:k3-malicious-4-t4-r2}
\end{subfigure}\hfill
\begin{subfigure}[t]{0.27\linewidth}
\centering
\includegraphics[width=\linewidth]{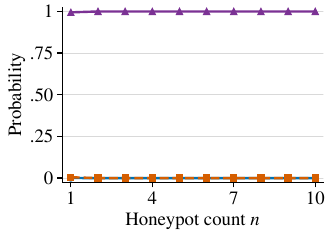}
\caption{$T=4$, $\lambda^m=2$}
\label{fig:k3-malicious-4-t4-r3}
\end{subfigure}\par\medskip
\begin{subfigure}[t]{0.27\linewidth}
\centering
\includegraphics[width=\linewidth]{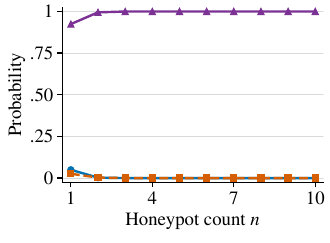}
\caption{$T=6$, $\lambda^m=0.5$}
\label{fig:k3-malicious-4-t6-r1}
\end{subfigure}\hfill
\begin{subfigure}[t]{0.27\linewidth}
\centering
\includegraphics[width=\linewidth]{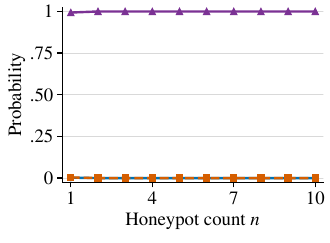}
\caption{$T=6$, $\lambda^m=1$}
\label{fig:k3-malicious-4-t6-r2}
\end{subfigure}\hfill
\begin{subfigure}[t]{0.27\linewidth}
\centering
\includegraphics[width=\linewidth]{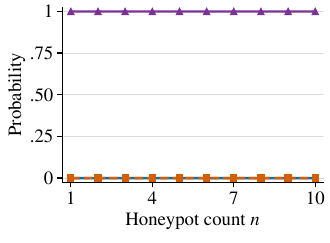}
\caption{$T=6$, $\lambda^m=2$}
\label{fig:k3-malicious-4-t6-r3}
\end{subfigure}\par\medskip
\caption{Transient state probabilities for malicious insiders with $K=3$, $\lambda^{\prime m}=2$, and $\lambda^{\prime\prime m}=4$.}
\label{fig: Prob_MC_malicious_K3_4}
\end{figure*}

Figures~\ref{fig: Prob_MC_benign_K2}--\ref{fig: Prob_MC_malicious_K3_4} report the state probabilities at time $T$ for a process starting from $S_1$. Each figure varies $n$ within a panel, $\lambda^{(\cdot)}$ across columns, and $T$ across rows. The three curves distinguish no recorded interactions, with probability $P^{(\cdot)}_1(T)$, suspicion without detection, with probability $P^{(\cdot)}_2(T)$ for $K=2$ and $P^{(\cdot)}_2(T)+P^{(\cdot)}_{3A}(T)+P^{(\cdot)}_{3B}(T)$ for $K=3$, and detection, with probability $P^{(\cdot)}_{K+1}(T)$. Detection represents a true detection for malicious insiders and a false detection for benign insiders. We evaluate the probabilities using the matrix exponential, including coincident-rate cases and the reduced chain for $K=3$, $n=1$. An independent calculation based on per-honeypot visit counts agrees with all 450 probability vectors to within $10^{-12}$.

Increasing $n$ or $T$ weakly increases detection probabilities for both insider types, subject to the zero-detection boundary cases identified in Section~\ref{sec:TheoreticalProperties}. The initial-state probability decreases as $e^{-n\lambda^{(\cdot)}T}$, but the probability of remaining suspicious can first increase and then decrease as interactions accumulate and detection becomes more likely. With $n$, $T$, and the interaction rates fixed, increasing the threshold from $K=2$ to $K=3$ leaves the initial-state probability unchanged and weakly reduces the detection probability. The probability removed from detection remains in the suspicious states. Thus, a larger suspicious-state probability under $K=3$ reflects the additional evidence required by the detection rule.

\begin{table}[tbp]
\caption{Configurations minimizing false detection subject to malicious-detection probability of at least $0.99$.}
\label{tab:balanceTn}
\centering
\small
\setlength{\tabcolsep}{4pt}
\begin{tabular}{cc ccc cc cc}
\toprule
 & & \multicolumn{3}{c}{Minimum honeypot count} & \multicolumn{2}{c}{Selected design} & \multicolumn{2}{c}{Detection probability} \\
\cmidrule(lr){3-5}\cmidrule(lr){6-7}\cmidrule(lr){8-9}
$\lambda^b$ & $\lambda^m$ & $n^*(2)$ & $n^*(4)$ & $n^*(6)$ & $n^*$ & $T^*$ & $P^m_{K+1}(T^*)$ & $P^b_{K+1}(T^*)$ \\
\midrule
\multicolumn{9}{l}{$K=2$} \\
\midrule
0.1 & 0.5 & 6 & 3 & 2 & 2 & 6 & 0.9959 & 0.2773 \\
0.1 & 1 & 3 & 2 & 1 & 1 & 6 & 0.9950 & 0.0672 \\
0.1 & 2 & 2 & 1 & 1 & 1 & 4 & 0.9970 & 0.0329 \\
0.2 & 0.5 & 6 & 3 & 2 & 2 & 6 & 0.9959 & 0.5562 \\
0.2 & 1 & 3 & 2 & 1 & 1 & 6 & 0.9950 & 0.1126 \\
0.2 & 2 & 2 & 1 & 1 & 1 & 4 & 0.9970 & 0.0581 \\
0.3 & 0.5 & 6 & 3 & 2 & 2 & 6 & 0.9959 & 0.7444 \\
0.3 & 1 & 3 & 2 & 1 & 1 & 6 & 0.9950 & 0.1441 \\
0.3 & 2 & 2 & 1 & 1 & 1 & 4 & 0.9970 & 0.0778 \\
\midrule
\multicolumn{9}{l}{$K=3$, $\lambda^{\prime\prime m}=0.5$} \\
\midrule
0.1 & 0.5 & 8 & 4 & 3 & 4 & 4 & 0.9960 & 0.1634 \\
0.1 & 1 & 4 & 3 & 2 & 4 & 2 & 0.9924 & 0.0335 \\
0.1 & 2 & 3 & 2 & 2 & 3 & 2 & 0.9991 & 0.0138 \\
0.2 & 0.5 & 8 & 4 & 3 & 4 & 4 & 0.9960 & 0.4655 \\
0.2 & 1 & 4 & 3 & 2 & 4 & 2 & 0.9924 & 0.1373 \\
0.2 & 2 & 3 & 2 & 2 & 3 & 2 & 0.9991 & 0.0578 \\
0.3 & 0.5 & 8 & 4 & 3 & 3 & 6 & 0.9984 & 0.6905 \\
0.3 & 1 & 4 & 3 & 2 & 2 & 6 & 0.9994 & 0.2198 \\
0.3 & 2 & 3 & 2 & 2 & 2 & 4 & 0.9998 & 0.1026 \\
\midrule
\multicolumn{9}{l}{$K=3$, $\lambda^{\prime\prime m}=4$} \\
\midrule
0.1 & 0.5 & 7 & 3 & 2 & 2 & 6 & 0.9946 & 0.0561 \\
0.1 & 1 & 4 & 2 & 1 & 1 & 6 & 0.9934 & 0 \\
0.1 & 2 & 3 & 1 & 1 & 1 & 4 & 0.9946 & 0 \\
0.2 & 0.5 & 7 & 3 & 2 & 2 & 6 & 0.9946 & 0.1447 \\
0.2 & 1 & 4 & 2 & 1 & 1 & 6 & 0.9934 & 0 \\
0.2 & 2 & 3 & 1 & 1 & 1 & 4 & 0.9946 & 0 \\
0.3 & 0.5 & 7 & 3 & 2 & 2 & 6 & 0.9946 & 0.2198 \\
0.3 & 1 & 4 & 2 & 1 & 1 & 6 & 0.9934 & 0 \\
0.3 & 2 & 3 & 1 & 1 & 1 & 4 & 0.9946 & 0 \\
\bottomrule
\end{tabular}
\par\smallskip
\begin{minipage}{0.98\linewidth}\footnotesize
All rows use $\lambda^{\prime b}=0.05$ and $\lambda^{\prime m}=2$. For $K=3$, $\lambda^{\prime\prime b}=0$. The selected count is $n^*=n^*(T^*)$. Probabilities are evaluated at $(n^*,T^*)$ and rounded only for display. Ties in false-detection probability are resolved by selecting the shorter horizon.
\end{minipage}
\end{table}

The effect of increasing the first-interaction rate depends on the baseline detection probability. For $K=2$, $n=4$, and $T=2$, an increase in $\lambda^b$ from $0.1$ to $0.2$ raises the false-detection probability from $0.1726$ to $0.4214$, a 2.44-fold increase [Figs.~\ref{fig:k2-benign-t2-r1} and~\ref{fig:k2-benign-t2-r2}]. Under the same deployment and horizon, increasing $\lambda^m$ from $0.5$ to $1$ raises true detection from $0.9585$ to $0.9985$ [Figs.~\ref{fig:k2-malicious-t2-r1} and~\ref{fig:k2-malicious-t2-r2}]. The absolute increase is smaller for malicious insiders, whose detection probability is already high, but their missed-detection probability falls from $0.0415$ to $0.0015$, approximately a 28-fold reduction. Near saturation, small changes in detection probability can therefore remain consequential when missed detections are costly. These comparisons also show why increasing interaction intensity indiscriminately can be undesirable when it increases benign engagement along with malicious engagement.

To examine the deployment implications of this saturation, we impose a common malicious-detection target of $0.99$. For each horizon, define
\[
n^*(T)=\min\big\{n\in\{1,\ldots,10\}:P^m_{K+1}(T)\ge0.99\big\}.
\]
At a fixed horizon, false-detection probability is nondecreasing in $n$, so the smallest count meeting the target also minimizes false detection. We then select $T^*$ from $\{2,4,6\}$ to minimize $P^b_{K+1}(T)$ evaluated at $n=n^*(T)$, and write $n^*=n^*(T^*)$. If horizons yield the same minimum false-detection probability, we select the shorter horizon. Table~\ref{tab:balanceTn} reports the minimum count for every horizon, the selected configuration, and its detection probabilities. This comparison isolates the trade-off between deployment size and monitoring duration before introducing the exposure and false-detection costs in Section~\ref{sec:theory_optimization}.

The target is attainable with at most eight honeypots at $T=2$ and at most three at $T=6$ in all considered scenarios. For example, when $\lambda^m=0.5$, the two-stage model requires six, three, and two honeypots at $T=2,4,6$, respectively. The three-stage model with $\lambda^{\prime\prime m}=0.5$ requires eight, four, and three. Raising $\lambda^{\prime\prime m}$ to $4$ reduces these counts to seven, three, and two. Faster third visits make the additional evidence requirement easier to satisfy, but the required count still depends on opportunities for first and second interactions. Since the target concerns only malicious detection, $n^*(T)$ is independent of the benign rates. Those rates affect the subsequent selection of $T^*$ through false-detection risk.

A longer horizon can reduce false detection under this selection rule by allowing a smaller deployment. It does not reduce false detection at a fixed $n$. The integer-valued count makes the trade-off uneven. For $K=2$, $\lambda^b=0.1$, and $\lambda^m=1$, extending $T$ from $2$ to $4$ reduces $n^*(T)$ from three to two, yet false detection rises from $0.1049$ to $0.1527$. Extending $T$ to $6$ permits one honeypot and lowers false detection to $0.0672$. Across the 27 scenarios, the selected horizon is $6$ in 14 cases, $4$ in nine, and $2$ in four. For two horizons requiring the same count, the longer horizon cannot improve false detection. For example, with $K=2$ and $\lambda^m=2$, one honeypot meets the target at both $T=4$ and $T=6$, and the shorter horizon is selected for every benign rate.

The three-stage model has a particularly useful boundary configuration. With one honeypot and $\lambda^{\prime\prime b}=0$, benign insiders cannot supply the third interaction, so false detection is exactly zero. When $\lambda^{\prime\prime m}=4$, one honeypot meets the malicious-detection target at $T=6$ for $\lambda^m=1$ and at $T=4$ for $\lambda^m=2$. For the latter rate, $T=6$ also yields zero false detection, so the stated tie rule selects $T=4$. These configurations exploit the difference in repeated engagement between insider types. Their zero false-detection probability depends on the assumption that benign third visits are impossible.

Requiring more evidence does not, however, guarantee a lower false-detection probability after deployment and horizon are adjusted to meet the same malicious-detection target. For $\lambda^b=0.3$ and $\lambda^m=1$, the selected false-detection probability is $0.1441$ under $K=2$, compared with $0.2198$ under $K=3$ when $\lambda^{\prime\prime m}=0.5$. The larger deployment needed by the three-stage model outweighs its reduction in false detection at a fixed configuration. Moreover, attaining the $0.99$ target alone does not ensure an acceptable false-detection rate. With $\lambda^b=0.3$ and $\lambda^m=0.5$, even the selected two-stage configuration falsely detects a benign insider with probability $0.7444$. These results motivate assessing detection timing and the expected costs of malicious exposure and false detections alongside the detection probabilities.

\subsection{Monotonicity and Sensitivity of Horizon-Conditioned Detection Time}
\label{sec:Experiment_theory_monotonicity_detection-time}

\begin{figure*}[tbp]
\centering
\captionsetup[sub]{font=small}
\includegraphics[width=0.4\linewidth]{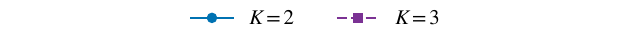}\par\smallskip
\begin{subfigure}[t]{0.27\linewidth}
\centering
\includegraphics[width=\linewidth]{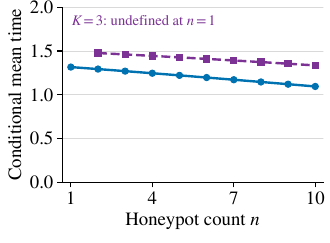}
\caption{$T=2$, $\lambda^b=0.1$}
\label{fig:time-benign-t2-r1}
\end{subfigure}\hfill
\begin{subfigure}[t]{0.27\linewidth}
\centering
\includegraphics[width=\linewidth]{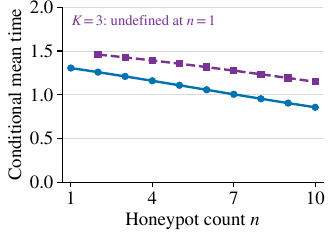}
\caption{$T=2$, $\lambda^b=0.2$}
\label{fig:time-benign-t2-r2}
\end{subfigure}\hfill
\begin{subfigure}[t]{0.27\linewidth}
\centering
\includegraphics[width=\linewidth]{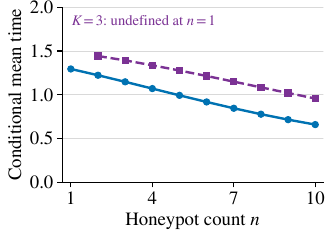}
\caption{$T=2$, $\lambda^b=0.3$}
\label{fig:time-benign-t2-r3}
\end{subfigure}\par\medskip
\begin{subfigure}[t]{0.27\linewidth}
\centering
\includegraphics[width=\linewidth]{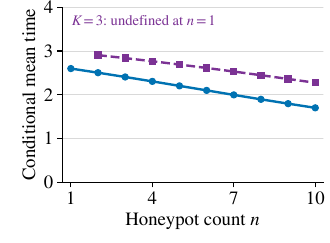}
\caption{$T=4$, $\lambda^b=0.1$}
\label{fig:time-benign-t4-r1}
\end{subfigure}\hfill
\begin{subfigure}[t]{0.27\linewidth}
\centering
\includegraphics[width=\linewidth]{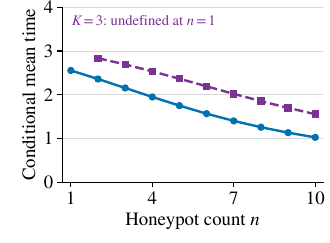}
\caption{$T=4$, $\lambda^b=0.2$}
\label{fig:time-benign-t4-r2}
\end{subfigure}\hfill
\begin{subfigure}[t]{0.27\linewidth}
\centering
\includegraphics[width=\linewidth]{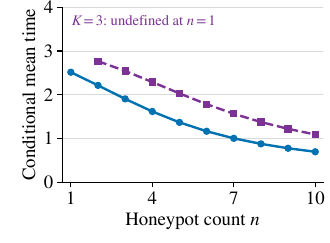}
\caption{$T=4$, $\lambda^b=0.3$}
\label{fig:time-benign-t4-r3}
\end{subfigure}\par\medskip
\begin{subfigure}[t]{0.27\linewidth}
\centering
\includegraphics[width=\linewidth]{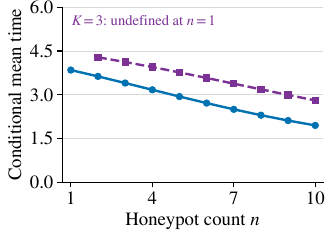}
\caption{$T=6$, $\lambda^b=0.1$}
\label{fig:time-benign-t6-r1}
\end{subfigure}\hfill
\begin{subfigure}[t]{0.27\linewidth}
\centering
\includegraphics[width=\linewidth]{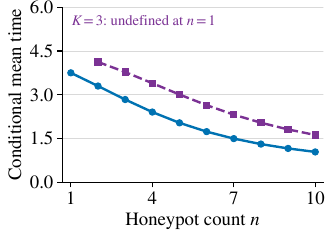}
\caption{$T=6$, $\lambda^b=0.2$}
\label{fig:time-benign-t6-r2}
\end{subfigure}\hfill
\begin{subfigure}[t]{0.27\linewidth}
\centering
\includegraphics[width=\linewidth]{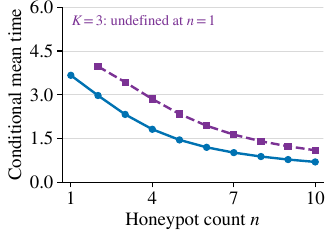}
\caption{$T=6$, $\lambda^b=0.3$}
\label{fig:time-benign-t6-r3}
\end{subfigure}\par\medskip
\caption{Expected detection time conditional on detection by $T$ for benign insiders, with $\lambda^{\prime b}=0.05$ and $\lambda^{\prime\prime b}=0$. For $K=3$, the conditional mean at $n=1$ is undefined and is omitted, since false detection is impossible.}
\label{fig:ExpectedTime_benign}
\end{figure*}

\begin{figure*}[tbp]
\centering
\captionsetup[sub]{font=small}
\includegraphics[width=0.4\linewidth]{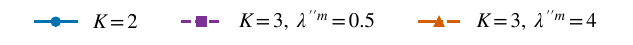}\par\smallskip
\begin{subfigure}[t]{0.27\linewidth}
\centering
\includegraphics[width=\linewidth]{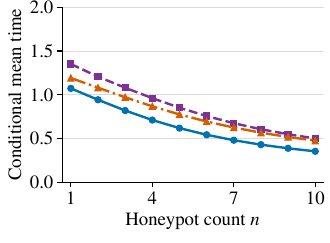}
\caption{$T=2$, $\lambda^m=0.5$}
\label{fig:time-malicious-t2-r1}
\end{subfigure}\hfill
\begin{subfigure}[t]{0.27\linewidth}
\centering
\includegraphics[width=\linewidth]{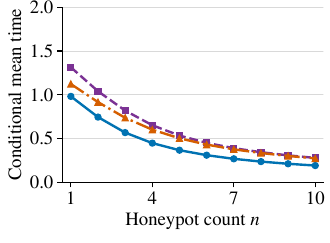}
\caption{$T=2$, $\lambda^m=1$}
\label{fig:time-malicious-t2-r2}
\end{subfigure}\hfill
\begin{subfigure}[t]{0.27\linewidth}
\centering
\includegraphics[width=\linewidth]{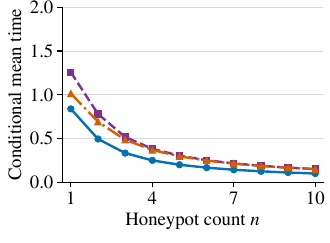}
\caption{$T=2$, $\lambda^m=2$}
\label{fig:time-malicious-t2-r3}
\end{subfigure}\par\medskip
\begin{subfigure}[t]{0.27\linewidth}
\centering
\includegraphics[width=\linewidth]{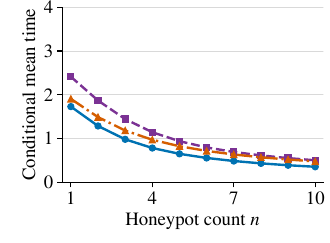}
\caption{$T=4$, $\lambda^m=0.5$}
\label{fig:time-malicious-t4-r1}
\end{subfigure}\hfill
\begin{subfigure}[t]{0.27\linewidth}
\centering
\includegraphics[width=\linewidth]{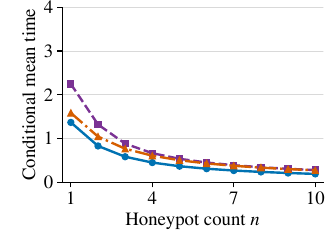}
\caption{$T=4$, $\lambda^m=1$}
\label{fig:time-malicious-t4-r2}
\end{subfigure}\hfill
\begin{subfigure}[t]{0.27\linewidth}
\centering
\includegraphics[width=\linewidth]{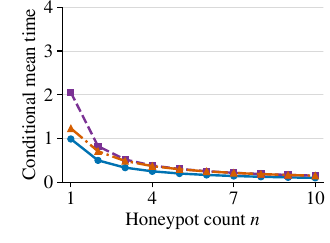}
\caption{$T=4$, $\lambda^m=2$}
\label{fig:time-malicious-t4-r3}
\end{subfigure}\par\medskip
\begin{subfigure}[t]{0.27\linewidth}
\centering
\includegraphics[width=\linewidth]{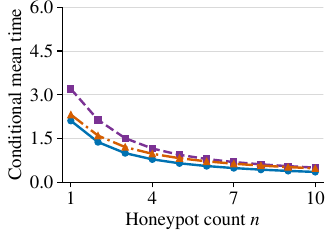}
\caption{$T=6$, $\lambda^m=0.5$}
\label{fig:time-malicious-t6-r1}
\end{subfigure}\hfill
\begin{subfigure}[t]{0.27\linewidth}
\centering
\includegraphics[width=\linewidth]{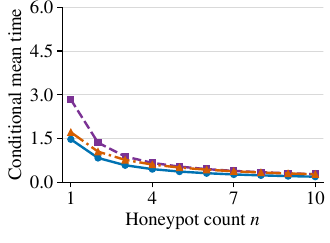}
\caption{$T=6$, $\lambda^m=1$}
\label{fig:time-malicious-t6-r2}
\end{subfigure}\hfill
\begin{subfigure}[t]{0.27\linewidth}
\centering
\includegraphics[width=\linewidth]{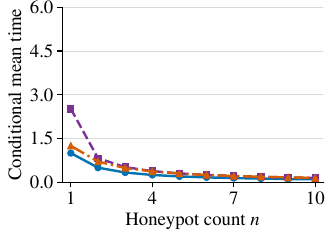}
\caption{$T=6$, $\lambda^m=2$}
\label{fig:time-malicious-t6-r3}
\end{subfigure}\par\medskip
\caption{Expected detection time conditional on detection by $T$ for malicious insiders, with $\lambda^{\prime m}=2$. %
}
\label{fig:ExpectedTime_malicious}
\end{figure*}

Figures~\ref{fig:ExpectedTime_benign} and~\ref{fig:ExpectedTime_malicious} report $\mathbb{E}[\tau^{(\cdot)}\mid\tau^{(\cdot)}\le T]$ for benign and malicious insiders, respectively. Each panel compares the detection thresholds at a common horizon and first-interaction rate. We evaluate the conditional means from the survival probabilities and their integrals, with an independent check obtained by integrating the detection-time density constructed from per-honeypot visit counts. For benign insiders under $K=3$, the value at $n=1$ is undefined because $\lambda^{\prime\prime b}=0$ makes detection impossible. The corresponding curves therefore begin at $n=2$.

Throughout the experimental grid, the three-stage model has a larger conditional mean than the two-stage model whenever both means are defined. Increasing $n$ or the first-interaction rate reduces the conditional mean for both thresholds. For malicious insiders under $K=3$, increasing $\lambda^{\prime\prime m}$ from $0.5$ to $4$ also reduces the mean at every fixed $(n,T,\lambda^m)$. These comparisons describe the timing of detections within the horizon. The conditional ordering across thresholds and the three-stage sensitivity patterns are numerical findings for the stated scenarios, since the pathwise detection-time comparisons in Section~\ref{sec:TheoreticalProperties} alone do not order means conditioned on different detection events.

The distinction matters because
\[
\mathbb{E}[\tau^{(\cdot)}\mid\tau^{(\cdot)}\le T]
=\frac{\mathbb{E}[\tau^{(\cdot)}\mathbf{1}\{\tau^{(\cdot)}\le T\}]}{\mathbb{P}(\tau^{(\cdot)}\le T)}.
\]
Requiring a third interaction changes both the accumulated detection times in the numerator and the probability in the denominator. A smaller detection probability alone cannot explain a larger conditional mean. The numerator is bounded above by $T\mathbb{P}(\tau^{(\cdot)}\le T)$, so every defined conditional mean remains between $0$ and $T$. For example, at $n=2$, $T=2$, and $\lambda^b=0.1$, the three-stage false-detection probability is only $0.0032$, while the conditional mean is $1.4770$ [Fig.~\ref{fig:time-benign-t2-r1}]. This mean describes the timing of a small subset of benign insiders. At larger deployments, false detections can instead be widespread. With $n=10$, $T=6$, and $\lambda^b=0.3$, the three-stage false-detection probability exceeds $0.9999$ and its conditional mean is $1.0964$ [Fig.~\ref{fig:time-benign-t6-r3}]. A short time to falsely flagging a benign insider therefore has a different operational meaning from rapid malicious detection. For benign insiders, the cost model prices the probability of false detection rather than its conditional timing.

The increase in conditional detection time with $T$ follows from the expanding set of detected trajectories. To make this effect explicit, fix the insider type, threshold, honeypot count, and rates, and write $F(T)=\mathbb{P}(\tau\le T)$, $f(T)=F'(T)$, and $\mu(T)=\mathbb{E}[\tau\mid\tau\le T]$. Differentiating the conditional first moment gives, whenever $F(T)>0$,
\begin{equation}
\frac{\mathrm{d}\mu(T)}{\mathrm{d}T}
=\frac{f(T)}{F(T)}\big(T-\mu(T)\big)\ge0.
\label{eq:conditional_mean_horizon_sensitivity}
\end{equation}
Additional detections enter the average at later times, while detections already occurring before the shorter horizon retain their original times. For $K=3$, $n=1$, $\lambda^m=0.5$, and $\lambda^{\prime\prime m}=0.5$, extending $T$ from $2$ to $6$ increases the conditional mean from $1.3514$ to $3.1985$, alongside an increase in detection probability from $0.1805$ to $0.7566$ [Figs.~\ref{fig:time-malicious-t2-r1} and~\ref{fig:time-malicious-t6-r1}]. The shorter-horizon mean represents the early detections among a population in which most malicious insiders remain undetected.

Horizon sensitivity becomes small when little detection probability remains beyond the shorter horizon. For $K=2$, $n=10$, and $\lambda^m=0.5$, the conditional mean changes from $0.3535$ at $T=2$ to $0.3538$ at $T=6$, close to the unconditional mean $0.3538$. Although saturation limits the benefit of further deployment for detection probability, additional honeypots can still shorten detection times. At $T=2$ and $\lambda^m=0.5$, increasing $n$ from six to ten raises the two-stage detection probability from $0.9928$ to $0.9998$ while reducing the conditional mean from $0.5427$ to $0.3535$. Thus, deployments that appear similar against a probability target can differ substantially in detection speed.

\begin{table}[tbp]
\caption{Detection timing at the minimum honeypot count meeting the $0.99$ malicious-detection target, for $\lambda^b=0.1$ and $\lambda^m=0.5$.}
\label{tab:detection_time_tradeoff}
\centering
\small
\begin{tabular}{ccccc}
\toprule
$n^*(T)$ & $T$ & $P^m_{K+1}(T)$ & $\mathbb{E}[\tau^m\mid\tau^m\le T]$ & $P^b_{K+1}(T)$ \\
\midrule
\multicolumn{5}{l}{$K=2$} \\
\midrule
6 & 2 & 0.9928 & 0.5427 & 0.3187 \\
3 & 4 & 0.9950 & 0.9818 & 0.2987 \\
2 & 6 & 0.9959 & 1.3768 & 0.2773 \\
\midrule
\multicolumn{5}{l}{$K=3$, $\lambda^{\prime\prime m}=0.5$} \\
\midrule
8 & 2 & 0.9954 & 0.6060 & 0.1904 \\
4 & 4 & 0.9960 & 1.1435 & 0.1634 \\
3 & 6 & 0.9984 & 1.5101 & 0.1842 \\
\midrule
\multicolumn{5}{l}{$K=3$, $\lambda^{\prime\prime m}=4$} \\
\midrule
7 & 2 & 0.9940 & 0.6252 & 0.1421 \\
3 & 4 & 0.9928 & 1.1824 & 0.0772 \\
2 & 6 & 0.9946 & 1.5989 & 0.0561 \\
\bottomrule
\end{tabular}
\par\smallskip
\begin{minipage}{0.94\linewidth}\footnotesize\centering
All rows use $\lambda^{\prime b}=0.05$ and $\lambda^{\prime m}=2$. For $K=3$, $\lambda^{\prime\prime b}=0$. %
\end{minipage}
\end{table}

Table~\ref{tab:detection_time_tradeoff} adds detection timing to the deployment--horizon comparison for $\lambda^b=0.1$ and $\lambda^m=0.5$. Every row uses the smallest count meeting the same $0.99$ malicious-detection target. Under $K=2$, extending $T$ from $2$ to $6$ permits a reduction from six honeypots to two and lowers false detection from $0.3187$ to $0.2773$. The malicious conditional mean, however, increases from $0.5427$ to $1.3768$. A longer horizon with fewer honeypots can therefore improve false-detection performance while increasing the time to malicious detection. This comparison changes the deployment as well as the horizon, in contrast to the fixed-deployment effect in \eqref{eq:conditional_mean_horizon_sensitivity}.

Faster repeated engagement can also be used to reduce deployment size rather than to improve detection speed. At $T=6$, the three-stage model requires three honeypots when $\lambda^{\prime\prime m}=0.5$ and two when $\lambda^{\prime\prime m}=4$. The corresponding conditional means are $1.5101$ and $1.5989$. Despite the faster third-visit rate, the smaller deployment yields a slightly longer conditional detection time, while reducing false detection from $0.1842$ to $0.0561$. These results show why the probability-based configurations in Table~\ref{tab:balanceTn} require a further cost assessment. The objective in \eqref{eq:defender_objective_time} combines the detection-weighted mean $P^m_{K+1}(T)\mathbb{E}[\tau^m\mid\tau^m\le T]$, exposure assigned to missed detections, and false-detection costs, so the conditional mean is one component of the design decision.

\subsection{Analysis of Cost-based Optimization}\label{sec:theory_optimization}

\begin{figure*}[tbp]
\centering
\captionsetup[sub]{font=small,justification=centering}
\includegraphics[width=0.4\linewidth]{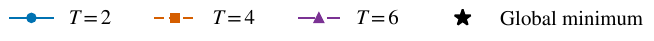}\par\smallskip
\begin{subfigure}[t]{0.27\linewidth}
\centering
\includegraphics[width=\linewidth]{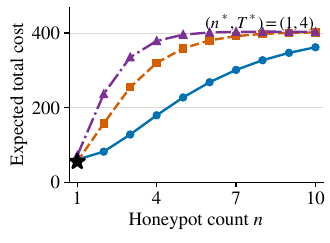}
\caption{$K=2$, $(c^a,c^f)=(10,100)$}
\label{fig:cost-model1-cost1}
\end{subfigure}\hfill
\begin{subfigure}[t]{0.27\linewidth}
\centering
\includegraphics[width=\linewidth]{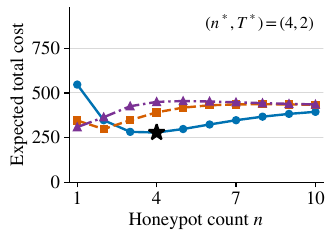}
\caption{$K=2$, $(c^a,c^f)=(100,100)$}
\label{fig:cost-model1-cost2}
\end{subfigure}\hfill
\begin{subfigure}[t]{0.27\linewidth}
\centering
\includegraphics[width=\linewidth]{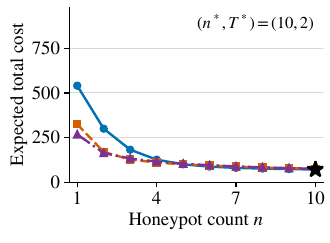}
\caption{$K=2$, $(c^a,c^f)=(100,10)$}
\label{fig:cost-model1-cost3}
\end{subfigure}\par\medskip
\begin{subfigure}[t]{0.27\linewidth}
\centering
\includegraphics[width=\linewidth]{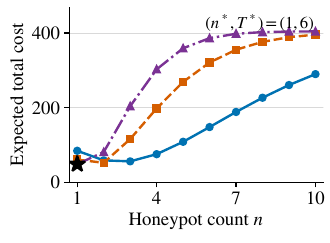}
\caption{$K=3$, $\lambda^{\prime\prime m}=0.5$, $(c^a,c^f)=(10,100)$}
\label{fig:cost-model2-cost1}
\end{subfigure}\hfill
\begin{subfigure}[t]{0.27\linewidth}
\centering
\includegraphics[width=\linewidth]{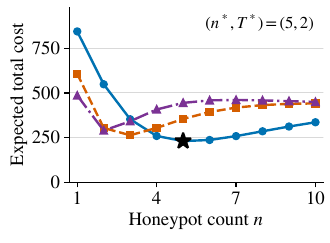}
\caption{$K=3$, $\lambda^{\prime\prime m}=0.5$, $(c^a,c^f)=(100,100)$}
\label{fig:cost-model2-cost2}
\end{subfigure}\hfill
\begin{subfigure}[t]{0.27\linewidth}
\centering
\includegraphics[width=\linewidth]{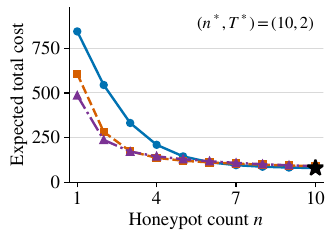}
\caption{$K=3$, $\lambda^{\prime\prime m}=0.5$, $(c^a,c^f)=(100,10)$}
\label{fig:cost-model2-cost3}
\end{subfigure}\par\medskip
\begin{subfigure}[t]{0.27\linewidth}
\centering
\includegraphics[width=\linewidth]{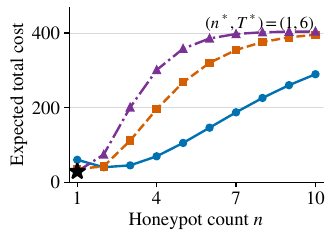}
\caption{$K=3$, $\lambda^{\prime\prime m}=4$, $(c^a,c^f)=(10,100)$}
\label{fig:cost-model3-cost1}
\end{subfigure}\hfill
\begin{subfigure}[t]{0.27\linewidth}
\centering
\includegraphics[width=\linewidth]{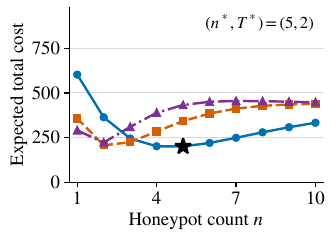}
\caption{$K=3$, $\lambda^{\prime\prime m}=4$, $(c^a,c^f)=(100,100)$}
\label{fig:cost-model3-cost2}
\end{subfigure}\hfill
\begin{subfigure}[t]{0.27\linewidth}
\centering
\includegraphics[width=\linewidth]{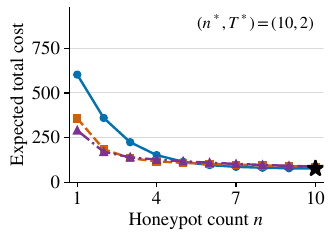}
\caption{$K=3$, $\lambda^{\prime\prime m}=4$, $(c^a,c^f)=(100,10)$}
\label{fig:cost-model3-cost3}
\end{subfigure}\par\medskip
\caption{Expected total cost for $\lambda^b=0.2$, $\lambda^m=0.5$, $\lambda^{\prime b}=0.05$, $\lambda^{\prime m}=2$, $Z=5$, and $\bar T=10$. For $K=3$, $\lambda^{\prime\prime b}=0$. Each star marks the minimum over all 30 feasible designs in its panel.}
\label{fig:cost_profiles}
\end{figure*}

We now evaluate the joint selection of $n$ and $T$ using the expected cost in \eqref{eq:defender_objective_time}. We fix $Z=5$ and $\bar T=10$ and consider $(c^a,c^f)=(10,100)$, $(100,100)$, and $(100,10)$. For each threshold and interaction-rate configuration, we enumerate all 30 designs with $n\in\{1,\ldots,10\}$ and $T\in\{2,4,6\}$. This yields 81 optimization scenarios. Unlike the probability-based comparisons in Tables~\ref{tab:balanceTn} and~\ref{tab:detection_time_tradeoff}, the cost model imposes no minimum malicious-detection probability. Tables~\ref{tab:balanceTn_K=2} and~\ref{tab:balanceTn_K=3} report the minimizing configurations and their expected costs. Figure~\ref{fig:cost_profiles} shows the underlying cost profiles for $\lambda^b=0.2$ and $\lambda^m=0.5$. Independent integration of the detection-time density confirms every minimizing design.

\begin{table}[tbp]
\caption{Cost-minimizing configurations and expected total costs for $K=2$.}
\label{tab:balanceTn_K=2}
\centering
\small
\setlength{\tabcolsep}{4pt}
\begin{tabular}{cc cc cc cc}
\toprule
 & & \multicolumn{2}{c}{$(c^a,c^f)=(10,100)$} & \multicolumn{2}{c}{$(c^a,c^f)=(100,100)$} & \multicolumn{2}{c}{$(c^a,c^f)=(100,10)$} \\
\cmidrule(lr){3-4}\cmidrule(lr){5-6}\cmidrule(lr){7-8}
$\lambda^b$ & $\lambda^m$ & $(n^*,T^*)$ & $C^*$ & $(n^*,T^*)$ & $C^*$ & $(n^*,T^*)$ & $C^*$ \\
\midrule
0.1 & 0.5 & $(1,4)$ & 45.36 & $(5,2)$ & 176.31 & $(10,2)$ & 58.73 \\
0.1 & 1 & $(1,4)$ & 29.96 & $(3,2)$ & 107.10 & $(8,2)$ & 41.95 \\
0.1 & 2 & $(1,2)$ & 20.40 & $(2,2)$ & 71.41 & $(6,2)$ & 29.41 \\
0.2 & 0.5 & $(1,4)$ & 55.47 & $(4,2)$ & 278.24 & $(10,2)$ & 71.38 \\
0.2 & 1 & $(1,2)$ & 39.37 & $(2,2)$ & 173.33 & $(10,2)$ & 54.94 \\
0.2 & 2 & $(1,2)$ & 23.57 & $(2,2)$ & 104.92 & $(4,2)$ & 41.86 \\
0.3 & 0.5 & $(1,4)$ & 63.32 & $(1,6)$ & 320.91 & $(10,2)$ & 74.66 \\
0.3 & 1 & $(1,2)$ & 42.16 & $(1,4)$ & 199.24 & $(10,2)$ & 58.22 \\
0.3 & 2 & $(1,2)$ & 26.36 & $(1,4)$ & 132.75 & $(10,2)$ & 49.13 \\
\bottomrule
\end{tabular}
\par\smallskip
\begin{minipage}{0.98\linewidth}\footnotesize \centering
All rows use $Z=5$, $\bar T=10$, $\lambda^{\prime b}=0.05$, and $\lambda^{\prime m}=2$.
\end{minipage}
\end{table}

Relative costs have a substantial effect on deployment size. Dividing the objective by $c^a$ expresses it as expected malicious exposure plus $4c^f/c^a$ times the benign false-detection probability. The three settings therefore assign weights of $40$, $4$, and $0.4$ to that probability, respectively. Even equal unit costs retain a fourfold multiplier because four employees are benign. Across this grid, increasing the relative cost of attack damage selects weakly larger deployments, lowers expected malicious exposure, and permits higher false-detection probabilities. For example, under $K=2$, $\lambda^b=0.1$, and $\lambda^m=0.5$, the selected design changes from $(1,4)$ to $(5,2)$ and then to $(10,2)$ across the three settings. The corresponding true-detection probabilities are $0.8197$, $0.9826$, and $0.9998$, while false detection rises from $0.0329$ to $0.2452$ and $0.5800$. The cost criterion thus balances missed and delayed malicious detection against disruption to benign employees, rather than enforcing near-certain detection in every scenario.

The shortest horizon, $T^*=2$, is selected in 65 of the 81 scenarios, with $T^*=4$ and $T^*=6$ each selected eight times. This pattern is particularly strong when attack damage is relatively expensive: all 27 scenarios with $(c^a,c^f)=(100,10)$ select $T^*=2$, compared with 15 of 27 when $(c^a,c^f)=(10,100)$. Equation~\eqref{eq:cost_horizon_derivative} clarifies that extending monitoring cannot increase expected malicious exposure at fixed $n$, but can raise total cost through additional false detections. For $K=2$, $n=1$, $\lambda^b=0.2$, $\lambda^m=0.5$, and $(c^a,c^f)=(10,100)$, total costs at $T=2,4,6$ are $60.77$, $55.47$, and $71.38$, respectively [Fig.~\ref{fig:cost-model1-cost1}]. An intermediate horizon is optimal even with the deployment held fixed.

\begin{table}[tbp]
\caption{Cost-minimizing configurations and expected total costs for $K=3$.}
\label{tab:balanceTn_K=3}
\centering
\small
\setlength{\tabcolsep}{4pt}
\begin{tabular}{cc cc cc cc}
\toprule
 & & \multicolumn{2}{c}{$(c^a,c^f)=(10,100)$} & \multicolumn{2}{c}{$(c^a,c^f)=(100,100)$} & \multicolumn{2}{c}{$(c^a,c^f)=(100,10)$} \\
\cmidrule(lr){3-4}\cmidrule(lr){5-6}\cmidrule(lr){7-8}
$\lambda^b$ & $\lambda^m$ & $(n^*,T^*)$ & $C^*$ & $(n^*,T^*)$ & $C^*$ & $(n^*,T^*)$ & $C^*$ \\
\midrule
\multicolumn{8}{l}{$\lambda^{\prime\prime m}=0.5$} \\
\midrule
0.1 & 0.5 & $(4,2)$ & 33.78 & $(7,2)$ & 134.27 & $(10,2)$ & 62.64 \\
0.1 & 1 & $(3,2)$ & 17.39 & $(5,2)$ & 79.50 & $(10,2)$ & 39.65 \\
0.1 & 2 & $(3,2)$ & 10.78 & $(4,2)$ & 51.63 & $(9,2)$ & 26.40 \\
0.2 & 0.5 & $(1,6)$ & 48.54 & $(5,2)$ & 230.36 & $(10,2)$ & 79.31 \\
0.2 & 1 & $(2,2)$ & 31.73 & $(4,2)$ & 127.21 & $(9,2)$ & 56.24 \\
0.2 & 2 & $(2,2)$ & 14.71 & $(3,2)$ & 75.67 & $(6,2)$ & 39.02 \\
0.3 & 0.5 & $(1,6)$ & 48.54 & $(5,2)$ & 312.96 & $(10,2)$ & 87.20 \\
0.3 & 1 & $(2,2)$ & 35.78 & $(3,2)$ & 169.20 & $(10,2)$ & 64.20 \\
0.3 & 2 & $(2,2)$ & 18.77 & $(3,2)$ & 103.15 & $(5,2)$ & 48.11 \\
\midrule
\multicolumn{8}{l}{$\lambda^{\prime\prime m}=4$} \\
\midrule
0.1 & 0.5 & $(2,4)$ & 26.49 & $(6,2)$ & 121.52 & $(10,2)$ & 59.85 \\
0.1 & 1 & $(3,2)$ & 14.60 & $(5,2)$ & 75.75 & $(10,2)$ & 39.13 \\
0.1 & 2 & $(2,2)$ & 9.17 & $(4,2)$ & 50.27 & $(8,2)$ & 26.28 \\
0.2 & 0.5 & $(1,6)$ & 29.02 & $(5,2)$ & 200.17 & $(10,2)$ & 76.52 \\
0.2 & 1 & $(1,6)$ & 17.70 & $(3,2)$ & 113.87 & $(8,2)$ & 55.43 \\
0.2 & 2 & $(2,2)$ & 12.25 & $(3,2)$ & 72.06 & $(6,2)$ & 38.62 \\
0.3 & 0.5 & $(1,6)$ & 29.02 & $(2,4)$ & 223.56 & $(10,2)$ & 84.41 \\
0.3 & 1 & $(1,6)$ & 17.70 & $(3,2)$ & 141.36 & $(10,2)$ & 63.69 \\
0.3 & 2 & $(1,6)$ & 12.51 & $(2,2)$ & 87.27 & $(5,2)$ & 47.42 \\
\bottomrule
\end{tabular}
\par\smallskip
\begin{minipage}{0.98\linewidth}\footnotesize \centering
All rows use $Z=5$, $\bar T=10$, $\lambda^{\prime b}=0.05$, and $\lambda^{\prime m}=2$. For $K=3$, $\lambda^{\prime\prime b}=0$.
\end{minipage}
\end{table}

The zero-false-detection boundary produces a different horizon choice. Under $K=3$, $n=1$, and $\lambda^{\prime\prime b}=0$, extending $T$ creates no additional false detections. It therefore reduces total cost by allowing more malicious detections. Whenever a single honeypot is selected under $K=3$, the optimizer chooses the longest feasible horizon, $T=6$. This also distinguishes the cost criterion from the tie rule in Table~\ref{tab:balanceTn}. For $\lambda^b=0.3$, $\lambda^m=2$, and $\lambda^{\prime\prime m}=4$, both $T=4$ and $T=6$ meet the $0.99$ target with one honeypot and zero false detection. Under $(c^a,c^f)=(10,100)$, however, extending $T$ lowers cost from $12.79$ to $12.51$, so the two horizons are economically distinct.

The higher evidence threshold reduces false detection at a fixed configuration but also delays malicious detection. After optimizing the deployment and horizon, $K=3$ has lower expected cost than $K=2$ in 22 of the 27 paired comparisons for each malicious third-visit rate. All five comparisons favoring $K=2$ occur under $(c^a,c^f)=(100,10)$, where the cost of slower malicious detection can outweigh the false-detection savings. For example, with $\lambda^b=0.2$ and $\lambda^m=0.5$, the optimized costs under $K=2$ and $K=3$ with $\lambda^{\prime\prime m}=4$ are $55.47$ and $29.02$ when $(c^a,c^f)=(10,100)$. Under $(100,10)$, the respective costs are $71.38$ and $76.52$ [Figs.~\ref{fig:cost-model1-cost3} and~\ref{fig:cost-model3-cost3}]. The preferred threshold therefore depends on the relative value of rapid malicious detection and protection against false detections.

Faster repeated malicious engagement improves the cost performance of the three-stage model. Increasing $\lambda^{\prime\prime m}$ from $0.5$ to $4$ leaves benign detection unchanged and reduces malicious exposure at every fixed design. The minimum achievable cost consequently cannot increase. The optimizer sometimes uses this improvement to reduce the honeypot count, as with $\lambda^b=0.1$, $\lambda^m=0.5$, and equal unit costs, where $n^*$ falls from seven to six. Other settings retain the same deployment and realize the benefit through faster detection. The joint decision can also produce a larger count when first interactions become faster. Under $K=3$, $\lambda^b=0.2$, $\lambda^{\prime\prime m}=0.5$, and $(c^a,c^f)=(10,100)$, increasing $\lambda^m$ from $0.5$ to $1$ changes the optimum from $(1,6)$ to $(2,2)$. The higher rate makes a shorter horizon with more honeypots preferable to a single honeypot monitored for longer.

Finally, 16 scenarios select the upper deployment limit $n^*=10$, all under $(c^a,c^f)=(100,10)$. Detection probabilities in these settings are already high, so the value of deployment also reflects shorter exposure before detection. The objective contains no separate deployment charge, making the feasible count limit consequential for these solutions. The reported configurations minimize cost over the stated grid, with every optimum unique at full numerical precision. Some neighboring designs are nearly equivalent in cost, however. For $K=3$, $\lambda^b=0.1$, $\lambda^m=2$, $\lambda^{\prime\prime m}=4$, and $(c^a,c^f)=(100,10)$, eight and nine honeypots at $T=2$ differ in expected cost by only $0.0079$. The tables identify the minimizing count, while the cost profiles show how strongly the objective distinguishes alternative configurations.

\subsection{Comparison with Detection at the First Interaction}\label{sec:K1_benchmark}

The preceding analysis compares thresholds that both require multiple interactions. We now introduce $K=1$ to assess whether the reduction in false detections justifies waiting beyond the first interaction. Appendix~\ref{app:K1_benchmark} derives its detection measures and expected cost.

Under $K=1$, both detection probabilities depend on $nT$. Any configuration meeting a malicious-detection target $P_2^m(T)\ge q$, with $0<q<1$, therefore satisfies
\begin{equation}
P_2^b(T)\ge 1-(1-q)^{\lambda^b/\lambda^m}.
\label{eq:K1_false_detection_bound}
\end{equation}
For $\lambda^b=0.1$, $\lambda^m=0.5$, and $q=0.99$, this bound is $0.6019$. On the stated grid, the minimum false-detection probability is $0.6321$ at $(n,T)=(5,2)$, compared with $0.2773$ for $K=2$ and $0.1634$ and $0.0561$ for $K=3$ with $\lambda^{\prime\prime m}=0.5$ and $4$, respectively (Table~\ref{tab:balanceTn}). The first-interaction rule cannot exploit differences in repeated engagement to reduce this bound.

For the cost comparison, each threshold selects $n$ and $T$ separately over the same 30 designs, without a minimum detection probability. Table~\ref{tab:K1_cost_comparison} summarizes cost reductions relative to optimized $K=1$. Both $K=2$ and $K=3$ have lower costs in all nine rate combinations under $(c^a,c^f)=(10,100)$ and $(100,100)$.

\begin{table}[htbp]
\caption{Median cost reduction (\%) relative to optimized $K=1$.}
\label{tab:K1_cost_comparison}
\centering
\small
\begin{tabular}{cccc}
\toprule
$(c^a,c^f)$ & $K=2$ & $K=3,\ \lambda^{\prime\prime m}=0.5$ & $K=3,\ \lambda^{\prime\prime m}=4$ \\
\midrule
$(10,100)$ & $74.0$ & $81.1$ & $88.3$ \\
$(100,100)$ & $42.0$ & $55.2$ & $59.9$ \\
$(100,10)$ & $-9.4$ & $-14.2$ & $-9.6$ \\
\bottomrule
\end{tabular}
\par\smallskip
\begin{minipage}{0.94\linewidth}\footnotesize
Each entry is the median of $100(C_1^*-C_K^*)/C_1^*$ over the nine $(\lambda^b,\lambda^m)$ pairs. Here $C_K^*$ is the minimum cost at threshold $K$. Negative values indicate higher costs than $K=1$.
\end{minipage}
\end{table}

Under $(100,10)$, however, $K=1$ has the lowest cost in six of nine combinations for either malicious third-interaction rate. Faster intervention then outweighs the additional false-detection costs. The remaining three combinations favor $K=3$. Thus, $K=2$ is never the least-cost threshold among all three on this grid, although it can outperform $K=3$ in the paired comparisons above.

\subsection{Managerial Implications}
\label{sec:PracticalInsights_new}

Detection probability alone is an incomplete basis for deployment decisions. Additional honeypots can shorten malicious exposure even when detection is already highly likely (Section~\ref{sec:Experiment_theory_monotonicity_detection-time}). Table~\ref{tab:detection_time_tradeoff} also shows that longer monitoring can meet the same $0.99$ detection target with fewer honeypots, but with a higher expected detection time among malicious insiders detected within the horizon. %

Requiring more evidence should be evaluated together with the deployment and monitoring duration needed to maintain detection performance. In Table~\ref{tab:balanceTn}, the false-detection advantage of requiring three interactions can reverse after configurations are adjusted to meet the same malicious-detection target. The cost comparisons in Tables~\ref{tab:balanceTn_K=2}, \ref{tab:balanceTn_K=3}, and~\ref{tab:K1_cost_comparison} show that the preferred threshold also depends on the relative costs of malicious exposure and falsely accusing benign insiders.

Monitoring duration has a different effect when deployment is held fixed. %
Proposition~\ref{prop:marginal_monitoring_value} shows that a small extension reduces expected cost when the savings from additional malicious detections exceed the additional false-detection costs across benign employees. The intermediate horizon selected in Fig.~\ref{fig:cost-model1-cost1} illustrates this balance. An increase in the mean time among detected insiders is therefore insufficient grounds for shortening monitoring.

\section{Conclusions}\label{sec:conclusions}

This study develops a stochastic framework for jointly choosing decoy deployment and monitoring duration when insider detection requires multiple interactions. The two- and three-stage models distinguish first visits from repeated visits to individual honeypots and evaluate benign and malicious insiders under the same detection rule. Closed-form detection probabilities and expected detection times conditional on detection within a finite horizon provide the performance measures. Boundary cases identify when restrictions on benign repeated visits preclude false detection, while rate-dominance conditions order malicious and benign detection times. The cost model connects these measures to decisions that balance the expected costs of malicious exposure and false detections. The numerical analysis shows that additional honeypots can reduce malicious exposure even when detection probability is already high. Comparing thresholds after adjusting deployment and monitoring also shows that requiring more evidence can reduce modeled expected cost, while high exposure costs relative to false-detection costs can favor detection at the first interaction.

Future research could extend the single-cycle framework to adaptive monitoring across successive cycles, using information from earlier cycles to update rate estimates and revise honeypot counts, placement, and types. Such models could jointly choose evidence thresholds and horizons, decide when to reset detection records, and determine how much evidence to retain or discount to balance persistent malicious activity against isolated benign interactions. Incorporating redeployment costs and game-theoretic models of malicious insider adaptation would support decisions over longer planning horizons. Extensions with multiple malicious insiders would also require assumptions about coordination and dependence across their actions.

\SingleSpacedXI
\footnotesize
\bibliographystyle{informs2014}
\bibliography{reference}
\normalsize
\SingleSpacedXI

\clearpage
\begin{APPENDICES}
\renewcommand{\theHsection}{appendix.\Alph{section}}
\renewcommand{\theHsubsection}{\theHsection.\arabic{subsection}}
\numberwithin{equation}{section}
\renewcommand{\theHequation}{\theHsection.\arabic{equation}}
\section{Analytical Derivations and Proofs}\label{app:analytical_results}

This appendix derives the detection measures used in Section~\ref{Sec:CTMC}, explains the cost formulation in Section~\ref{Sec:Optimization}, and proves the results in Section~\ref{sec:TheoreticalProperties}. Throughout, first-interaction rates are positive, repeated-interaction rates are nonnegative, and honeypot counts are positive integers. Conditional detection times are considered only for $T>0$ and positive detection probability.

\subsection{Conditional Expected Detection Time for the Two-stage Model}
\label{App:K2_CondMeanDerivation}

We derive \eqref{eq:K2_conditional_expected_detection_time} for a fixed insider type $(\cdot)\in\{b,m\}$ and a horizon $T>0$. We suppress the type superscript on aggregate rates and probabilities. First assume that both stage rates are positive. Boundary cases are addressed at the end of this subsection.

\paragraph{Holding-time representation.}

For $K=2$, the CTMC has the unidirectional structure $S_1\rightarrow S_2\rightarrow S_3$, where $S_3$ is absorbing. Starting from $S_1$, the absorption time $\tau$ can be written as
\begin{equation}
\label{eq:tau_sum_two_exponentials}
\tau = X_1 + X_2,
\end{equation}
where $X_1$ is the sojourn time in $S_1$ and $X_2$ is the sojourn time in $S_2$. By the CTMC
assumption, $X_1$ and $X_2$ are independent exponential random variables with rates
\begin{equation}
\label{eq:k2_rates_def}
X_1 \sim \mathrm{Exp}(\lambda_1), \qquad
X_2 \sim \mathrm{Exp}(\lambda_2),
\end{equation}
where
\begin{equation}
\label{eq:k2_l1_l2}
\lambda_1 = n\lambda^{(\cdot)},\qquad
\lambda_2 = (n-1)\lambda^{(\cdot)} + \lambda^{\prime(\cdot)}.
\end{equation}

Because the detection process is evaluated over a finite horizon $T$, the expected detection time
\emph{among detected insiders} is
\begin{equation}
\label{eq:cond_mean_target}
\mathbb{E}[\tau \mid \tau \le T].
\end{equation}

\paragraph{Conditioning on detection within the horizon.}

Let $\overline F(t):=\mathbb{P}(\tau>t)$ denote the survival function of $\tau$. Using the identity
$\min(\tau,T)=\int_0^T \mathbf{1}\{\tau>t\}\,dt$ and Tonelli's theorem, we have
\begin{equation}
\label{eq:min_identity}
\mathbb{E}[\min(\tau,T)]
=
\int_0^T \mathbb{P}(\tau>t)\,dt
=
\int_0^T \overline F(t)\,dt.
\end{equation}
Also,
\begin{equation}
\label{eq:min_decomp}
\mathbb{E}[\min(\tau,T)]
=
\mathbb{E}\!\left[\tau\,\mathbf{1}\{\tau\le T\}\right]
+
T\,\mathbb{P}(\tau>T)
=
\mathbb{E}\!\left[\tau\,\mathbf{1}\{\tau\le T\}\right]
+
T\,\overline F(T).
\end{equation}
Rearranging \eqref{eq:min_decomp} gives
\begin{equation}
\label{eq:tau_indicator}
\mathbb{E}\!\left[\tau\,\mathbf{1}\{\tau\le T\}\right]
=
\mathbb{E}[\min(\tau,T)] - T\,\overline F(T).
\end{equation}
Finally, dividing by $\mathbb{P}(\tau\le T)=1-\overline F(T)$ yields
\begin{equation}
\label{eq:cond_mean_in_terms_survival}
\mathbb{E}[\tau \mid \tau \le T]
=
\frac{\mathbb{E}[\min(\tau,T)] - T\,\overline F(T)}{1 - \overline F(T)}
=
\frac{\displaystyle \int_0^T \overline F(t)\,dt - T\,\overline F(T)}{1 - \overline F(T)}.
\end{equation}
Thus, it suffices to obtain closed forms for $\overline F(t)$ and $\int_0^T \overline F(t)\,dt$.

\paragraph{Survival function.}

For $t\ge 0$, conditioning on $X_1$ and using independence in \eqref{eq:tau_sum_two_exponentials}
gives
\begin{align}
\overline F(t)
=
\mathbb{P}(X_1+X_2>t)
&=
\int_0^\infty \mathbb{P}(X_2>t-x)\,f_{X_1}(x)\,dx \nonumber\\
&=
\int_0^t \mathbb{P}(X_2>t-x)\,\lambda_1 e^{-\lambda_1 x}\,dx
\;+\;
\int_t^\infty 1\cdot \lambda_1 e^{-\lambda_1 x}\,dx,
\label{eq:survival_split}
\end{align}
where $f_{X_1}(x)=\lambda_1 e^{-\lambda_1 x}$.

For $0\le x\le t$, we have $\mathbb{P}(X_2>t-x)=e^{-\lambda_2(t-x)}$, so the first integral becomes
\begin{align}
\int_0^t e^{-\lambda_2(t-x)}\,\lambda_1 e^{-\lambda_1 x}\,dx
&=
\lambda_1 e^{-\lambda_2 t}\int_0^t e^{(\lambda_2-\lambda_1)x}\,dx \nonumber\\
&=
\lambda_1 e^{-\lambda_2 t}\cdot
\frac{e^{(\lambda_2-\lambda_1)t}-1}{\lambda_2-\lambda_1}
=
\frac{\lambda_1}{\lambda_2-\lambda_1}\left(e^{-\lambda_1 t}-e^{-\lambda_2 t}\right),
\label{eq:int1_eval}
\end{align}
where $\lambda_1\neq\lambda_2$. Equal positive rates are handled below.

The second integral in \eqref{eq:survival_split} is
\begin{equation}
\label{eq:int2_eval}
\int_t^\infty \lambda_1 e^{-\lambda_1 x}\,dx = e^{-\lambda_1 t}.
\end{equation}
Combining \eqref{eq:int1_eval} and \eqref{eq:int2_eval} yields
\begin{align}
\overline F(t)
&=
\frac{\lambda_1}{\lambda_2-\lambda_1}\left(e^{-\lambda_1 t}-e^{-\lambda_2 t}\right)
+
e^{-\lambda_1 t} \nonumber\\
&=
\frac{\lambda_2 e^{-\lambda_1 t}-\lambda_1 e^{-\lambda_2 t}}{\lambda_2-\lambda_1}.
\label{eq:survival_closed_form}
\end{align}

\paragraph{Integrated survival and the conditional mean.}

Using \eqref{eq:survival_closed_form}, we integrate term-by-term to obtain
\begin{align}
\int_0^T \overline F(t)\,dt
&=
\frac{1}{\lambda_2-\lambda_1}
\left[
\lambda_2\int_0^T e^{-\lambda_1 t}\,dt
-
\lambda_1\int_0^T e^{-\lambda_2 t}\,dt
\right] \nonumber\\
&=
\frac{1}{\lambda_2-\lambda_1}
\left[
\lambda_2\cdot \frac{1-e^{-\lambda_1 T}}{\lambda_1}
-
\lambda_1\cdot \frac{1-e^{-\lambda_2 T}}{\lambda_2}
\right].
\label{eq:int_survival_closed_form}
\end{align}
Moreover, evaluating \eqref{eq:survival_closed_form} at $T$ gives
\begin{equation}
\label{eq:survival_at_T}
\overline F(T)=\frac{\lambda_2 e^{-\lambda_1 T}-\lambda_1 e^{-\lambda_2 T}}{\lambda_2-\lambda_1}.
\end{equation}

Substituting \eqref{eq:int_survival_closed_form} and \eqref{eq:survival_at_T} into
\eqref{eq:cond_mean_in_terms_survival} yields
\begin{equation}
\label{eq:K2_final_conditional_mean}
\mathbb{E}[\tau \mid \tau \le T]
=
\frac{
\displaystyle
\frac{1}{\lambda_2-\lambda_1}
\left[
\frac{\lambda_2}{\lambda_1}\left(1-e^{-\lambda_1 T}\right)
-
\frac{\lambda_1}{\lambda_2}\left(1-e^{-\lambda_2 T}\right)
\right]
-
T\,
\frac{\lambda_2 e^{-\lambda_1 T}-\lambda_1 e^{-\lambda_2 T}}{\lambda_2-\lambda_1}
}{
\displaystyle
1-
\frac{\lambda_2 e^{-\lambda_1 T}-\lambda_1 e^{-\lambda_2 T}}{\lambda_2-\lambda_1}
}.
\end{equation}

Finally, reinstating the insider-type notation using \eqref{eq:k2_l1_l2} gives the expression in
\eqref{eq:K2_final_conditional_mean} with
$\lambda_1=\lambda^{(\cdot)}_1=n\lambda^{(\cdot)}$ and
$\lambda_2=\lambda^{(\cdot)}_2=(n-1)\lambda^{(\cdot)}+\lambda^{\prime(\cdot)}$.

\paragraph{Coincident and zero rates.}
If $\lambda_1=\lambda_2=r>0$, the absorption time has an Erlang distribution with two stages. Its survival function is $e^{-rt}(1+rt)$, which gives
\begin{equation}
\mathbb{E}[\tau\mid\tau\le T]
=\frac{2}{r}\,
\frac{1-e^{-rT}\big(1+rT+(rT)^2/2\big)}
{1-e^{-rT}(1+rT)}.
\label{eq:app_K2_equal_rates}
\end{equation}
This is also the limit of the complete expression in \eqref{eq:K2_final_conditional_mean}. If $n=1$ and $\lambda^{\prime(\cdot)}=0$, then $\lambda_2=0$ and the chain remains in $S_2$ after its first interaction. Detection is impossible, and the conditional mean is undefined. The positive-rate exponential representation above is not used in this case. More generally, the survival identity \eqref{eq:cond_mean_in_terms_survival} applies whenever detection by $T$ has positive probability, including distributions with mass at infinity. In that setting, $\tau\mathbf{1}\{\tau\le T\}$ denotes the variable equal to $\tau$ on $\{\tau\le T\}$ and zero otherwise.

\subsection{Transient Probabilities for the Three-stage Model}\label{app:K3_transient_probabilities}

For $n\ge2$, consider the CTMC with ordered state space $(S_1,S_2,S_{3A},S_{3B},S_4)$ and initial condition
\[
P^{(\cdot)}_1(0)=1, \qquad 
P^{(\cdot)}_2(0)=P^{(\cdot)}_{3A}(0)=P^{(\cdot)}_{3B}(0)=P^{(\cdot)}_4(0)=0.
\]
For insider type $(\cdot)\in\{b,m\}$, define the transition rates
\[
\lambda^{(\cdot)}_1 := n\lambda^{(\cdot)}, \qquad
\lambda^{(\cdot)}_2(A):=(n-1)\lambda^{(\cdot)}, \qquad
\lambda^{(\cdot)}_2(B):=\lambda^{\prime(\cdot)},
\]
\[
\lambda^{(\cdot)}_2:=\lambda^{(\cdot)}_2(A)+\lambda^{(\cdot)}_2(B),
\]
\[
\lambda^{(\cdot)}_3(A):=(n-2)\lambda^{(\cdot)}+2\lambda^{\prime(\cdot)}, \qquad
\lambda^{(\cdot)}_3(B):=(n-1)\lambda^{(\cdot)}+\lambda^{\prime\prime(\cdot)}.
\]

The transient probabilities satisfy the Kolmogorov forward equations
\begin{align}
\frac{\mathrm{d}}{\mathrm{d}t}P^{(\cdot)}_1(t)
&=
-\lambda^{(\cdot)}_1 P^{(\cdot)}_1(t), \label{eq:K3_app_P1}\\
\frac{\mathrm{d}}{\mathrm{d}t}P^{(\cdot)}_2(t)
&=
\lambda^{(\cdot)}_1 P^{(\cdot)}_1(t)
-\lambda^{(\cdot)}_2 P^{(\cdot)}_2(t), \label{eq:K3_app_P2}\\
\frac{\mathrm{d}}{\mathrm{d}t}P^{(\cdot)}_{3A}(t)
&=
\lambda^{(\cdot)}_2(A) P^{(\cdot)}_2(t)
-\lambda^{(\cdot)}_3(A) P^{(\cdot)}_{3A}(t), \label{eq:K3_app_P3A}\\
\frac{\mathrm{d}}{\mathrm{d}t}P^{(\cdot)}_{3B}(t)
&=
\lambda^{(\cdot)}_2(B) P^{(\cdot)}_2(t)
-\lambda^{(\cdot)}_3(B) P^{(\cdot)}_{3B}(t), \label{eq:K3_app_P3B}\\
\frac{\mathrm{d}}{\mathrm{d}t}P^{(\cdot)}_4(t)
&=
\lambda^{(\cdot)}_3(A) P^{(\cdot)}_{3A}(t)
+\lambda^{(\cdot)}_3(B) P^{(\cdot)}_{3B}(t). \label{eq:K3_app_P4}
\end{align}

Solving \eqref{eq:K3_app_P1} yields
\begin{equation}
P^{(\cdot)}_1(t)=e^{-\lambda^{(\cdot)}_1 t}.
\label{eq:K3_app_P1_sol}
\end{equation}
Substituting \eqref{eq:K3_app_P1_sol} into \eqref{eq:K3_app_P2} and using an integrating factor gives, for $\lambda^{(\cdot)}_1\ne\lambda^{(\cdot)}_2$,
\begin{equation}
P^{(\cdot)}_2(t)=
\frac{\lambda^{(\cdot)}_1}{\lambda^{(\cdot)}_2-\lambda^{(\cdot)}_1}
\left(
e^{-\lambda^{(\cdot)}_1 t}-e^{-\lambda^{(\cdot)}_2 t}
\right).
\label{eq:K3_app_P2_sol}
\end{equation}

Equations \eqref{eq:K3_app_P3A} and \eqref{eq:K3_app_P3B} admit the integral representations
\begin{align}
P^{(\cdot)}_{3A}(t)
&=
e^{-\lambda^{(\cdot)}_3(A)t}
\int_0^t
\lambda^{(\cdot)}_2(A) P^{(\cdot)}_2(s)
e^{\lambda^{(\cdot)}_3(A)s}\,\mathrm{d}s,
\label{eq:K3_app_P3A_int}\\
P^{(\cdot)}_{3B}(t)
&=
e^{-\lambda^{(\cdot)}_3(B)t}
\int_0^t
\lambda^{(\cdot)}_2(B) P^{(\cdot)}_2(s)
e^{\lambda^{(\cdot)}_3(B)s}\,\mathrm{d}s.
\label{eq:K3_app_P3B_int}
\end{align}
Substituting \eqref{eq:K3_app_P2_sol} and evaluating the integrals yields the following expressions when their denominators are nonzero:
\begin{align}
P^{(\cdot)}_{3A}(t)
&=
\frac{\lambda^{(\cdot)}_2(A)\lambda^{(\cdot)}_1}{\lambda^{(\cdot)}_2-\lambda^{(\cdot)}_1}
\left[
\frac{e^{-\lambda^{(\cdot)}_3(A)t}-e^{-\lambda^{(\cdot)}_1 t}}
{\lambda^{(\cdot)}_1-\lambda^{(\cdot)}_3(A)}
-
\frac{e^{-\lambda^{(\cdot)}_3(A)t}-e^{-\lambda^{(\cdot)}_2 t}}
{\lambda^{(\cdot)}_2-\lambda^{(\cdot)}_3(A)}
\right], \label{eq:K3_app_P3A_sol}\\
P^{(\cdot)}_{3B}(t)
&=
\frac{\lambda^{(\cdot)}_2(B)\lambda^{(\cdot)}_1}{\lambda^{(\cdot)}_2-\lambda^{(\cdot)}_1}
\left[
\frac{e^{-\lambda^{(\cdot)}_3(B)t}-e^{-\lambda^{(\cdot)}_1 t}}
{\lambda^{(\cdot)}_1-\lambda^{(\cdot)}_3(B)}
-
\frac{e^{-\lambda^{(\cdot)}_3(B)t}-e^{-\lambda^{(\cdot)}_2 t}}
{\lambda^{(\cdot)}_2-\lambda^{(\cdot)}_3(B)}
\right]. \label{eq:K3_app_P3B_sol}
\end{align}

Finally, by probability conservation,
\begin{equation}
P^{(\cdot)}_4(t)
=
1-
P^{(\cdot)}_1(t)
-
P^{(\cdot)}_2(t)
-
P^{(\cdot)}_{3A}(t)
-
P^{(\cdot)}_{3B}(t).
\label{eq:K3_app_P4_sol}
\end{equation}

These equations yield \eqref{eq:K3_transient_probabilities}. If $\lambda^{(\cdot)}_1=\lambda^{(\cdot)}_2=r>0$, then $P^{(\cdot)}_2(t)=rt e^{-rt}$. For any further coincident rates, the integral representations \eqref{eq:K3_app_P3A_int}--\eqref{eq:K3_app_P3B_int} remain valid. They can be evaluated directly or recovered as limits of the complete probability expressions. A branch with zero entry rate has zero probability, and a branch with zero exit rate retains any probability entering it. These cases also follow directly from the forward equations or the matrix exponential.

When $n=1$, Case A is impossible. Omit $S_{3A}$ and its incident transitions before evaluating any rates. The remaining chain is $S_1\to S_2\to S_{3B}\to S_4$, with rates $\lambda^{(\cdot)}$, $\lambda^{\prime(\cdot)}$, and $\lambda^{\prime\prime(\cdot)}$. The same forward-equation argument applies to this reduced chain, with $P^{(\cdot)}_{3A}(t)=0$. The expression $(n-2)\lambda^{(\cdot)}+2\lambda^{\prime(\cdot)}$ is not evaluated for the omitted state.

\subsection{Conditional Expected Detection Time for the Three-stage Model}
\label{app:K3_conditional_detection_time}

This appendix derives the closed-form expressions used in
\eqref{eq:K3_survival_mixture}--\eqref{eq:K3_conditional_expected_detection_time}
for the three-stage model ($K=3$). We first consider $n\ge2$, a positive total departure rate from $S_2$, and positive exit rates on all branches with positive probability. The remaining cases are treated below.

\paragraph{Branch decomposition.}
Starting from $S_1$, the process must transition to $S_2$ and then to either $S_{3A}$ or $S_{3B}$.
Let $X^{(\cdot)}_1$ denote the sojourn time in $S_1$, and let $X^{(\cdot)}_2$ denote the sojourn time in $S_2$.
By the CTMC construction,
\[
X^{(\cdot)}_1 \sim \mathrm{Exp}\!\big(\lambda^{(\cdot)}_1\big),
\qquad
\lambda^{(\cdot)}_1=n\lambda^{(\cdot)},
\]
and, since the total departure rate from $S_2$ is
\[
\lambda^{(\cdot)}_2=\lambda^{(\cdot)}_2(A)+\lambda^{(\cdot)}_2(B)
=(n-1)\lambda^{(\cdot)}+\lambda^{\prime(\cdot)},
\]
we have
\[
X^{(\cdot)}_2 \sim \mathrm{Exp}\!\big(\lambda^{(\cdot)}_2\big).
\]
Moreover, conditional on leaving $S_2$, the probability of branching into Case $A$ (new honeypot) or Case $B$ (repeat honeypot) is
\[
p_A^{(\cdot)}=\frac{\lambda^{(\cdot)}_2(A)}{\lambda^{(\cdot)}_2}
=\frac{(n-1)\lambda^{(\cdot)}}{(n-1)\lambda^{(\cdot)}+\lambda^{\prime(\cdot)}},
\qquad
p_B^{(\cdot)}=\frac{\lambda^{(\cdot)}_2(B)}{\lambda^{(\cdot)}_2}
=\frac{\lambda^{\prime(\cdot)}}{(n-1)\lambda^{(\cdot)}+\lambda^{\prime(\cdot)}}.
\]
Let $J\in\{A,B\}$ denote the realized branch, with $\mathbb{P}(J=j)=p_j^{(\cdot)}$. The joint density of the holding time and the selected branch is
\[
\mathbb{P}(X^{(\cdot)}_2\in\mathrm{d}x,J=j)
=\lambda^{(\cdot)}_2(j)e^{-\lambda^{(\cdot)}_2x}\,\mathrm{d}x
=p_j^{(\cdot)}\lambda^{(\cdot)}_2e^{-\lambda^{(\cdot)}_2x}\,\mathrm{d}x.
\]
Thus, the branch choice is independent of the holding time in $S_2$. It is also independent of the holding time in $S_1$ by the CTMC construction. This independence permits conditioning on the branch without changing either of the first two holding-time distributions.

\paragraph{Holding times conditional on a branch.}
Given a branch $J=j\in\{A,B\}$ with $p_j^{(\cdot)}>0$, the final transition to the absorbing state $S_4$ occurs after an additional exponential holding time
\[
X^{(\cdot)}_{3j}\sim \mathrm{Exp}\!\big(\lambda^{(\cdot)}_3(j)\big),
\]
where
\[
\lambda^{(\cdot)}_3(A)=(n-2)\lambda^{(\cdot)}+2\lambda^{\prime(\cdot)},
\qquad
\lambda^{(\cdot)}_3(B)=(n-1)\lambda^{(\cdot)}+\lambda^{\prime\prime(\cdot)}.
\]
By the CTMC construction and the branch independence established above, the holding times
$X^{(\cdot)}_1$, $X^{(\cdot)}_2$, and $X^{(\cdot)}_{3j}$ are independent for each fixed $j$.
Thus, the absorption time satisfies
\[
\tau^{(\cdot)}\mid (J=j)\overset{d}{=}X^{(\cdot)}_1+X^{(\cdot)}_2+X^{(\cdot)}_{3j},
\qquad j\in\{A,B\},
\]
which is a three-phase hypoexponential distribution.

\paragraph{Conditional survival function.}
Fix $j\in\{A,B\}$ and write
\[
a:=\lambda^{(\cdot)}_1,\quad b:=\lambda^{(\cdot)}_2,\quad c:=\lambda^{(\cdot)}_3(j).
\]
Assume $a,b,c$ are positive and pairwise distinct. Coincident rates are addressed below.
Let $\overline F^{(\cdot)}_j(t)=\mathbb{P}(\tau^{(\cdot)}>t\mid J=j)$.
The Laplace transform of the absorption-time density is $abc/[(s+a)(s+b)(s+c)]$. Partial-fraction inversion followed by integration of the density over $(t,\infty)$ yields
\begin{equation}
\label{eq:app_K3_survival_conditional}
\overline F^{(\cdot)}_j(t)
=
c^{(\cdot)}_{1j}e^{-a t}+c^{(\cdot)}_{2j}e^{-b t}+c^{(\cdot)}_{3j}e^{-c t},
\end{equation}
where the coefficients are
\begin{equation}
\label{eq:app_K3_coeffs}
c^{(\cdot)}_{1j}
=
\frac{bc}{(b-a)(c-a)},
\qquad
c^{(\cdot)}_{2j}
=
\frac{ac}{(a-b)(c-b)},
\qquad
c^{(\cdot)}_{3j}
=
\frac{ab}{(a-c)(b-c)}.
\end{equation}
\paragraph{Mixture over visit histories.}
By the law of total probability,
\begin{equation}
\label{eq:app_K3_survival_mixture}
\mathbb{P}\!\left(\tau^{(\cdot)} > t\right)
=
\sum_{j\in\{A,B\}} \mathbb{P}(J=j)\,\mathbb{P}\!\left(\tau^{(\cdot)} > t \mid J=j\right)
=
p_A^{(\cdot)}\,\overline F_A^{(\cdot)}(t)+p_B^{(\cdot)}\,\overline F_B^{(\cdot)}(t),
\end{equation}
which is \eqref{eq:K3_survival_mixture} together with \eqref{eq:K3_survival_components}.

\paragraph{Conditional first moment.}
Applying the survival identity established in Appendix~\ref{App:K2_CondMeanDerivation} gives
\begin{equation}
\label{eq:app_trunc_identity}
\mathbb{E}\!\left[\tau^{(\cdot)}\mathbf{1}\{\tau^{(\cdot)}\le T\}\right]
=\int_0^T\mathbb{P}(\tau^{(\cdot)}>t)\,\mathrm{d}t
-T\mathbb{P}(\tau^{(\cdot)}>T).
\end{equation}
This subtraction removes the contribution of trajectories that remain undetected at $T$ from the mean censored time.

\paragraph{Integrated survival.}
From \eqref{eq:app_K3_survival_mixture} and the exponential form in
\eqref{eq:app_K3_survival_conditional},
\[
\int_0^T \mathbb{P}\!\left(\tau^{(\cdot)} > t\right)\,dt
=
\sum_{j\in\{A,B\}} p_j^{(\cdot)} \int_0^T \overline F^{(\cdot)}_j(t)\,dt
=
\sum_{j\in\{A,B\}} p_j^{(\cdot)}
\left[
\frac{c^{(\cdot)}_{1j}}{a}\big(1-e^{-aT}\big)
+
\frac{c^{(\cdot)}_{2j}}{b}\big(1-e^{-bT}\big)
+
\frac{c^{(\cdot)}_{3j}}{c}\big(1-e^{-cT}\big)
\right].
\]
Substituting back $a=\lambda^{(\cdot)}_1$, $b=\lambda^{(\cdot)}_2$, and $c=\lambda^{(\cdot)}_3(j)$
yields \eqref{eq:K3_integrated_survival}.

\paragraph{Expected time among detections.}
Finally,
\[
\mathbb{E}\!\left[\tau^{(\cdot)} \mid \tau^{(\cdot)} \le T\right]
=
\frac{\mathbb{E}\!\left[\tau^{(\cdot)}\mathbf{1}\{\tau^{(\cdot)}\le T\}\right]}
{\mathbb{P}\!\left(\tau^{(\cdot)}\le T\right)}
=
\frac{
\displaystyle
\int_0^T \mathbb{P}\!\left(\tau^{(\cdot)} > t\right)\,dt
-
T\,\mathbb{P}\!\left(\tau^{(\cdot)} > T\right)
}{
\displaystyle
1-\mathbb{P}\!\left(\tau^{(\cdot)} > T\right)
}.
\]
Using $P^{(\cdot)}_4(T)=\mathbb{P}(\tau^{(\cdot)}\le T)$ and
$\mathbb{P}(\tau^{(\cdot)} > T)=1-P^{(\cdot)}_4(T)$ gives
\eqref{eq:K3_conditional_expected_detection_time}.

\paragraph{Coincident rates, unreachable branches, and nonabsorption.}
If rates coincide within a branch, take limits of the complete survival and integrated-survival expressions. Individual coefficients in \eqref{eq:app_K3_coeffs} can diverge and should not be assigned separate limits. Equality of the final rates on different branches does not cause a singularity unless one of those rates also equals the first or second stage rate. When all three rates within a branch equal $r>0$, its survival function is $e^{-rt}(1+rt+(rt)^2/2)$.

Branches with zero probability are omitted, so no conditional distribution is assigned to an impossible branch. If $\lambda^{(\cdot)}_2=0$, the chain cannot leave $S_2$ and detection is impossible. If a branch with positive probability has zero final-stage rate, set its conditional survival function to one. Its contribution to integrated survival is then $p_j^{(\cdot)}T$, and its contribution to the truncated first moment is zero. The conditional mean is defined only when $P^{(\cdot)}_4(T)>0$.

For $n=1$, omit Case A and use the reduced chain from Appendix~\ref{app:K3_transient_probabilities}. If $\lambda^{\prime(\cdot)}>0$, the only branch has probability one and successive rates $\lambda^{(\cdot)}$, $\lambda^{\prime(\cdot)}$, and $\lambda^{\prime\prime(\cdot)}$. Detection is impossible if either repeated-interaction rate is zero.

If the first two stage rates and every positive-probability branch's final rate are positive, absorption occurs with probability one and the unconditional mean is finite. Independence of the holding times and averaging over branches give
\[
\mathbb{E}[\tau^{(\cdot)}]
=\frac{1}{\lambda^{(\cdot)}_1}
+\frac{1}{\lambda^{(\cdot)}_2}
+\sum_{\substack{j\in\{A,B\}\\p_j^{(\cdot)}>0}}
\frac{p_j^{(\cdot)}}{\lambda^{(\cdot)}_3(j)}.
\]
The conditional mean converges to this value as $T\to\infty$, since the truncated first moment converges to $\mathbb{E}[\tau^{(\cdot)}]$ and the detection probability converges to one.

\subsection{Proof of Proposition~\ref{prop:marginal_monitoring_value}}
\label{app:cost_derivation}

Fix $n$, $K$, the interaction-rate parameters, and the exposure cutoff
$\bar T$. Write $F_m(t)=P^m_{K+1}(t)$ and
$F_b(t)=P^b_{K+1}(t)$. Under the exposure convention in
Section~\ref{Sec:Optimization}, the malicious insider's assigned
exposure duration is
\begin{equation}
L_T=
\begin{cases}
\tau^m, & \tau^m\le T,\\
\bar T, & \tau^m>T.
\end{cases}
\label{eq:app_exposure_duration}
\end{equation}
Hence,
\begin{equation}
\mathbb{E}[L_T]
=\int_0^T t\,\mathrm{d}F_m(t)
+\bar T\big(1-F_m(T)\big).
\label{eq:app_expected_exposure}
\end{equation}
For positive detection probability, the integral equals
$F_m(T)\mathbb{E}[\tau^m\mid\tau^m\le T]$.
For zero detection probability, it is zero. This representation
therefore covers the boundary cases without assigning a value
to an undefined conditional mean.

For each benign employee, the probability of a false detection
by $T$ is $F_b(T)$. Linearity of expectation gives an expected
number of false detections equal to $(Z-1)F_b(T)$, regardless
of dependence between employees. Multiplying
\eqref{eq:app_expected_exposure} by $c^a$ and adding
$(Z-1)c^fF_b(T)$ gives the objective in
\eqref{eq:defender_objective_time}.

To examine the effect of the horizon, let
$0<T_1<T_2\le\bar T$. Subtracting the two exposure expectations
gives
\begin{equation}
\mathbb{E}[L_{T_2}]-\mathbb{E}[L_{T_1}]
=\int_{T_1}^{T_2}(t-\bar T)\,\mathrm{d}F_m(t)
\le0.
\label{eq:app_exposure_horizon_order}
\end{equation}
An additional detection at time $t\in(T_1,T_2]$ replaces the
assigned missed-detection exposure $\bar T$ by $t$, while
detections before $T_1$ retain their original times. The
false-detection cost changes by
$(Z-1)c^f[F_b(T_2)-F_b(T_1)]\ge0$.
Thus, extending the horizon cannot increase expected malicious
exposure, although it can increase total cost through additional
false detections. The conditional mean detection time may
increase because the longer horizon includes later detections,
which does not contradict \eqref{eq:app_exposure_horizon_order}.

The CTMC detection probabilities are differentiable for $t>0$,
with derivatives $f_m(t)$ and $f_b(t)$. Holding $\bar T$ fixed,
differentiation of \eqref{eq:app_expected_exposure} yields
\[
\frac{\mathrm{d}\mathbb{E}[L_T]}{\mathrm{d}T}
=T f_m(T)-\bar T f_m(T)
=-(\bar T-T)f_m(T).
\]
Multiplying by $c^a$ and adding the derivative of the
false-detection cost gives \eqref{eq:cost_horizon_derivative}.
Its sign is negative precisely when
$c^a(\bar T-T)f_m(T)>(Z-1)c^f f_b(T)$.
Finally, at an interior minimizer over a continuous feasible
interval, differentiability requires the cost derivative to
vanish. The marginal saving and marginal false-detection cost
must therefore be equal, establishing the stated necessary
condition.\hfill$\square$

\subsection{Proof of Lemma~\ref{lem:n1_no_false_detection}}\label{app:proof_no_false_detection}

When $n=1$ and $\lambda^{\prime b}=0$, the second transition rate is
\[
\lambda^b_2=(n-1)\lambda^b+\lambda^{\prime b}=0.
\]
After the first interaction, the benign process remains in $S_2$. The detection state $S_3$ is therefore unreachable, giving $P^b_3(T)=0$ for every $T\ge0$.\hfill$\square$

\subsection{Proof of Proposition~\ref{Prop:DetectionTime}}\label{app:proof_K2_order}

\paragraph{Coupling and detection probabilities.}
For each insider type, the two aggregate rates are
\[
\lambda^{(\cdot)}_1=n\lambda^{(\cdot)},\qquad
\lambda^{(\cdot)}_2=(n-1)\lambda^{(\cdot)}+\lambda^{\prime(\cdot)}.
\]
The assumptions imply $\lambda_i^m\ge\lambda_i^b$ for $i=1,2$. Let $E_1,E_2$ be independent exponential variables with rate one. For a positive stage rate, define $X_i^{(\cdot)}=E_i/\lambda_i^{(\cdot)}$, and set $X_i^{(\cdot)}=\infty$ when the rate is zero. Each process then has the correct holding-time distributions, with independence between its two stages. Sharing $E_i$ across types gives
\[
X_i^m\le X_i^b\quad(i=1,2),\qquad
\tau^m=X_1^m+X_2^m\le X_1^b+X_2^b=\tau^b
\quad\text{almost surely under this coupling}.
\]
Consequently, $\{\tau^b\le T\}\subseteq\{\tau^m\le T\}$ under the coupling. Since coupling preserves each marginal distribution, it follows that $\tau^m\le_{\mathrm{st}}\tau^b$ and $P^m_3(T)\ge P^b_3(T)$ for every $T\ge0$.

\paragraph{Detection times conditional on the horizon.}
Now assume $\lambda_2^b>0$ and $T>0$, so all four stage rates are positive and both conditioning events have positive probability. The preceding pathwise comparison alone is insufficient because the two conditioning events differ. We instead compare the densities of the sums of holding times.

Let $g_{a,b}(t)$ be the density of the sum of independent exponential variables with positive rates $a$ and $b$. For $t>0$, convolution gives
\begin{equation}
g_{a,b}(t)=ab e^{-bt}\int_0^t e^{-(a-b)x}\,\mathrm{d}x.
\label{eq:app_K2_convolution_density}
\end{equation}
Fix $b>0$ and let $a_2\ge a_1>0$. Then
\begin{equation}
\frac{g_{a_2,b}(t)}{g_{a_1,b}(t)}
=\frac{a_2}{a_1}
\frac{\int_0^t e^{-(a_2-a_1)x}e^{-(a_1-b)x}\,\mathrm{d}x}
{\int_0^t e^{-(a_1-b)x}\,\mathrm{d}x}.
\label{eq:app_K2_density_ratio}
\end{equation}
The fraction is a weighted average of the nonincreasing function $e^{-(a_2-a_1)x}$ over $[0,t]$, with positive weights $e^{-(a_1-b)x}$. Increasing $t$ adds values no larger than the existing average, so the ratio is nonincreasing in $t$. By symmetry, increasing $b$ with $a$ fixed has the same effect. Increasing the two rates successively therefore shows that $g_m(t)/g_b(t)$ is nonincreasing, where $g_m$ and $g_b$ are the malicious and benign detection-time densities. This argument also covers coincident positive rates.

For $0<s<T$ and $0<u\le s<v<T$, the density-ratio ordering gives
\[
g_m(u)g_b(v)\ge g_b(u)g_m(v).
\]
Integrating over $u\in(0,s)$ and $v\in(s,T)$, and writing $F_m$ and $F_b$ for the two detection-time distribution functions, yields
\[
F_m(s)\big(F_b(T)-F_b(s)\big)
\ge F_b(s)\big(F_m(T)-F_m(s)\big).
\]
After canceling the common product and dividing by $F_m(T)F_b(T)>0$,
\begin{equation}
\frac{F_m(s)}{F_m(T)}\ge\frac{F_b(s)}{F_b(T)}.
\label{eq:app_K2_truncated_distribution_order}
\end{equation}
The two sides are the distribution functions conditional on detection by $T$. The inequality also holds at the endpoints, so these conditional detection times are stochastically ordered. Integrating their survival functions over $[0,T]$ gives
\[
\mathbb{E}[\tau^m\mid\tau^m\le T]
\le\mathbb{E}[\tau^b\mid\tau^b\le T],
\]
which proves \eqref{eq:K2_conditional_order}.\hfill$\square$

\subsection{Proof of Proposition~\ref{Prop:DetectionTimeK3} and the Three-stage Boundary Cases}\label{app:proof_K3_order}

Fix $n\ge1$ and consider the underlying interaction mechanism across the honeypots. We construct both insider types on a common probability space while preserving their marginal interaction processes.
For insider type $(\cdot)\in\{b,m\}$ and honeypot $i\in\{1,\dots,n\}$, let $C_i^{(\cdot)}(t)\in\{0,1,2,3\}$
denote the number of interactions with honeypot $i$ by time $t$, capped at three.
Define the total number of honeypot interactions by time $t$ as
\[
C^{(\cdot)}(t):=\sum_{i=1}^n C_i^{(\cdot)}(t),
\]
and define the absorption (detection) time as the time of the third interaction:
\[
\tau^{(\cdot)}:=\inf\{t\ge 0:\ C^{(\cdot)}(t)\ge 3\}.
\]
Under the $K=3$ detection rule, $\tau^{(\cdot)}$ coincides with the absorption time of the CTMC
on $\{S_1,S_2,S_{3A},S_{3B},S_4\}$ for $n\ge2$, or the reduced chain for $n=1$. Counts are capped at three per honeypot, since further visits cannot affect the time of the third total interaction. The auxiliary counts may continue after detection, but the detection time itself is unchanged.

\medskip
\paragraph{Shared Poisson clocks.}
Choose a constant $\nu>0$ such that
\[
\nu \ge \max\big\{\lambda^{m},\lambda^{\prime m},\lambda^{\prime\prime m},
                 \lambda^{b},\lambda^{\prime b},\lambda^{\prime\prime b}\big\}.
\]
For each honeypot $i$, let $\{N_i(t)\}_{t\ge 0}$ be an independent Poisson process with rate $\nu$,
with jump times $0<t_{i,1}<t_{i,2}<\cdots$, and let $\{U_{i,k}\}_{k\ge 1}$ be independent $\mathrm{Unif}(0,1)$
random variables, independent across $i,k$ and independent of all Poisson processes.

At each candidate time $t_{i,k}$, define the per-honeypot acceptance rates for insider type $(\cdot)$ as a function of
the current count $C_i^{(\cdot)}(t_{i,k}^-)$:
\[
r_{c}^{(\cdot)}:=
\begin{cases}
\lambda^{(\cdot)}, & c=0,\\
\lambda^{\prime(\cdot)}, & c=1,\\
\lambda^{\prime\prime(\cdot)}, & c=2,\\
0, & c\ge 3.
\end{cases}
\]
We declare that an actual interaction of type $(\cdot)$ occurs at $t_{i,k}$ if
\[
U_{i,k}\le \frac{r_{\,C_i^{(\cdot)}(t_{i,k}^-)}^{(\cdot)}}{\nu},
\]
in which case we set $C_i^{(\cdot)}(t_{i,k})=C_i^{(\cdot)}(t_{i,k}^-)+1$. Otherwise, we keep
$C_i^{(\cdot)}(t_{i,k})=C_i^{(\cdot)}(t_{i,k}^-)$.
Between candidate times, the processes $C_i^{(\cdot)}(\cdot)$ remain constant.

For each type, thinning the rate-$\nu$ clock with acceptance probability $r_c^{(\cdot)}/\nu$ gives an interaction rate $r_c^{(\cdot)}$ whenever the count is $c$. Independent clocks and marks across honeypots preserve the assumed independence of their interaction processes. Sharing the clocks and marks across insider types provides the coupling needed for the comparison.

\medskip
\paragraph{Ordering the interaction counts.}
Under the rate ordering \eqref{eq:rate_assumptions_K3_prop}, we have
\[
\lambda^{m}\ge\lambda^{b},\qquad
\lambda^{\prime m}\ge \lambda^{\prime b},\qquad
\lambda^{\prime\prime m}\ge \lambda^{\prime\prime b}.
\]
Hence, for each $c\in\{0,1,2,3\}$,
\[
r_c^{m}\ge r_c^{b} \qquad \Longrightarrow \qquad \frac{r_c^{m}}{\nu}\ge \frac{r_c^{b}}{\nu}.
\]
We claim that for each honeypot $i$,
\begin{equation}\label{eq:honeypot_domination_K3}
C_i^{m}(t)\ \ge\ C_i^{b}(t)\qquad \text{for all }t\ge 0,\ \text{almost surely}.
\end{equation}
The proof is by induction over the ordered candidate times $\{t_{i,k}\}_{k\ge 1}$.
At time $0$, both counts are $0$. Suppose \eqref{eq:honeypot_domination_K3} holds just before $t_{i,k}$.
If $C_i^{m}(t_{i,k}^-)>C_i^{b}(t_{i,k}^-)$, then after the update at $t_{i,k}$ the inequality still holds since
each count can increase by at most $1$.
If $C_i^{m}(t_{i,k}^-)=C_i^{b}(t_{i,k}^-)=c$, then acceptance for the benign process implies
$U_{i,k}\le r_c^{b}/\nu\le r_c^{m}/\nu$, and therefore acceptance for the malicious process as well.
Thus the inequality is preserved at $t_{i,k}$ in all cases, proving \eqref{eq:honeypot_domination_K3}.

Summing \eqref{eq:honeypot_domination_K3} over $i=1,\dots,n$ yields the total-count domination
\begin{equation}\label{eq:total_count_domination_K3}
C^{m}(t)=\sum_{i=1}^n C_i^{m}(t)\ \ge\ \sum_{i=1}^n C_i^{b}(t)=C^{b}(t)
\qquad \text{for all }t\ge 0,\ \text{almost surely}.
\end{equation}

\medskip
\paragraph{Detection-time and probability orderings.}
By definition,
\[
\tau^{(\cdot)}=\inf\{t\ge 0:\ C^{(\cdot)}(t)\ge 3\}.
\]
Therefore, \eqref{eq:total_count_domination_K3} implies the pathwise ordering
\begin{equation}\label{eq:tau_pathwise_order_K3}
\tau^{m}\le \tau^{b}\qquad \text{almost surely under the coupling}.
\end{equation}
In particular, $\tau^{m}\le_{\mathrm{st}}\tau^{b}$ (first-order stochastic dominance).
For every $T\ge0$, \eqref{eq:tau_pathwise_order_K3} implies the event inclusion
\[
\{\tau^{b}\le T\}\subseteq \{\tau^{m}\le T\}\qquad \text{almost surely},
\]
and taking probabilities gives
\[
P_4^{m}(T)=\mathbb{P}(\tau^{m}\le T)\ \ge\ \mathbb{P}(\tau^{b}\le T)=P_4^{b}(T),
\]
which proves \eqref{eq:P4_ordering_K3_prop}.
\hfill$\square$

\paragraph{Boundary cases for benign detection.}
For a finite horizon $T>0$, a finite sequence of transitions with positive rates has positive probability of being completed by $T$. Thus, it suffices to determine whether a path to three interactions exists. If $\lambda^{\prime b}=0$, each honeypot can contribute at most one interaction, so absorption is possible exactly when $n\ge3$. If $\lambda^{\prime b}>0$ and $\lambda^{\prime\prime b}=0$, each honeypot can contribute at most two interactions, so absorption is possible exactly when $n\ge2$. If both repeated-interaction rates are positive, a single honeypot can supply all three interactions. These observations prove the boundary conditions in Section~\ref{sec:theory_K3}. In particular, when $\lambda^{\prime b}=0$, the value of $\lambda^{\prime\prime b}$ is irrelevant because the second visit cannot occur.

\section{First-Interaction Detection Benchmark}\label{app:K1_benchmark}

For $K=1$, the first honeypot interaction triggers detection. For insider type $(\cdot)\in\{b,m\}$, the first-interaction times at the $n$ honeypots are independent exponential random variables with rate $\lambda^{(\cdot)}>0$. Their minimum has rate $r^{(\cdot)}=n\lambda^{(\cdot)}$. The CTMC therefore consists of $S_1\rightarrow S_2$, where $S_2$ is absorbing, with
\begin{equation}
P_1^{(\cdot)}(T)=e^{-r^{(\cdot)}T},
\qquad
P_2^{(\cdot)}(T)=1-e^{-r^{(\cdot)}T}.
\label{eq:K1_probabilities}
\end{equation}
Repeated-visit rates do not enter these expressions because detection precedes any repeated visit.

Integrating the density $r^{(\cdot)}e^{-r^{(\cdot)}t}$ gives
\begin{equation}
\int_0^T t\,r^{(\cdot)}e^{-r^{(\cdot)}t}\,\mathrm{d}t
=\frac{1-e^{-r^{(\cdot)}T}}{r^{(\cdot)}}-Te^{-r^{(\cdot)}T}.
\label{eq:K1_truncated_moment}
\end{equation}
Dividing by $P_2^{(\cdot)}(T)>0$ yields
\begin{equation}
\mathbb{E}[\tau^{(\cdot)}\mid\tau^{(\cdot)}\le T]
=\frac{1}{r^{(\cdot)}}-\frac{T}{e^{r^{(\cdot)}T}-1}.
\label{eq:K1_conditional_mean}
\end{equation}
Substitution into \eqref{eq:defender_objective_time}, with the same fixed exposure cutoff $\bar T\ge T$, gives
\begin{equation}
C_1(n,T)
=c^a\left[\frac{1-e^{-r^mT}}{r^m}+(\bar T-T)e^{-r^mT}\right]
+(Z-1)c^f\left(1-e^{-r^bT}\right).
\label{eq:K1_cost}
\end{equation}

Eliminating $nT$ from \eqref{eq:K1_probabilities} gives
\begin{equation}
P_2^b(T)=1-\bigl(1-P_2^m(T)\bigr)^{\lambda^b/\lambda^m}.
\label{eq:K1_probability_relation}
\end{equation}
The right-hand side is increasing in $P_2^m(T)$, proving \eqref{eq:K1_false_detection_bound}. In particular, equal first-interaction rates give identical true- and false-detection probabilities. For feasible configurations with a common product $nT$, \eqref{eq:K1_conditional_mean} also shows that the conditional mean is inversely proportional to $n$. Such configurations have identical detection probabilities, but larger deployments reduce the time to detection and the associated exposure costs when $c^a>0$.

Table~\ref{tab:K1_optima} reports all 27 cost-minimizing configurations under the experimental settings in Sections~\ref{Sec:ExprimentSetup} and~\ref{sec:theory_optimization}. The closed-form calculations and independent numerical integration of the detection-time density give the same minimizing design in every case.

\begin{table}[htbp]
\caption{Cost-minimizing honeypot counts and expected costs for $K=1$.}
\label{tab:K1_optima}
\centering
\small
\setlength{\tabcolsep}{5pt}
\begin{tabular}{cc cc cc cc}
\toprule
 & & \multicolumn{2}{c}{$(c^a,c^f)=(10,100)$} & \multicolumn{2}{c}{$(c^a,c^f)=(100,100)$} & \multicolumn{2}{c}{$(c^a,c^f)=(100,10)$} \\
\cmidrule(lr){3-4}\cmidrule(lr){5-6}\cmidrule(lr){7-8}
$\lambda^b$ & $\lambda^m$ & $n^*$ & $C_1^*$ & $n^*$ & $C_1^*$ & $n^*$ & $C_1^*$ \\
\midrule
$0.1$ & $0.5$ & $1$ & $114.58$ & $3$ & $283.65$ & $10$ & $54.62$ \\
$0.1$ & $1$ & $1$ & $91.98$ & $2$ & $195.61$ & $7$ & $44.42$ \\
$0.1$ & $2$ & $1$ & $78.88$ & $1$ & $136.24$ & $4$ & $34.53$ \\
$0.2$ & $0.5$ & $1$ & $173.94$ & $3$ & $382.70$ & $10$ & $59.30$ \\
$0.2$ & $1$ & $1$ & $151.35$ & $2$ & $284.01$ & $10$ & $49.27$ \\
$0.2$ & $2$ & $1$ & $138.25$ & $1$ & $195.61$ & $10$ & $44.27$ \\
$0.3$ & $0.5$ & $1$ & $222.55$ & $10$ & $419.04$ & $10$ & $59.94$ \\
$0.3$ & $1$ & $1$ & $199.95$ & $2$ & $343.26$ & $10$ & $49.90$ \\
$0.3$ & $2$ & $1$ & $186.85$ & $1$ & $244.21$ & $10$ & $44.90$ \\
\bottomrule
\end{tabular}
\par\smallskip
\begin{minipage}{0.94\linewidth}\footnotesize
All optima use $T^*=2$. The search uses $n\in\{1,\ldots,10\}$, $T\in\{2,4,6\}$, $Z=5$, and $\bar T=10$, with no minimum detection probability or separate deployment charge.
\end{minipage}
\end{table}

\end{APPENDICES}

\SingleSpacedXI

\end{document}